\documentclass[a4paper,10pt]{elsarticle}
\biboptions{sort&compress} %
\newcommand{\pdiff}[2]{\ensuremath{\frac{\partial #1}{\partial #2}}} 

\usepackage[english]{babel}
\usepackage[utf8]{inputenc}
\usepackage{graphicx}
\usepackage{caption}
\usepackage[labelformat=simple]{subcaption}

\usepackage{mathtools}
\usepackage{amsfonts}
\usepackage{xfrac}
\usepackage{siunitx}
\usepackage{placeins}
\usepackage{tikz}
\usepackage{pgfplots}
\usetikzlibrary{calc}
\usetikzlibrary{backgrounds}
\usepgflibrary{shapes}
\usetikzlibrary{through}
\usetikzlibrary{intersections}
\usetikzlibrary{matrix}
\usepackage{multirow}
\usepackage{nameref}
\usepackage{hhline}
\usetikzlibrary {arrows.meta}

\usepackage{hyperref}
\usepackage{cleveref}
\crefname{equation}{eq.}{eqs.}
\crefname{figure}{fig.}{figs.}
\Crefname{figure}{Figure}{Figures}

\crefname{table}{Table}{Tables}
\Crefname{table}{Table}{Tables}

\usepackage{xifthen}
\usepackage{wrapfig}

\usepackage{esdiff}
\usepackage{ifthen}
\usepackage{algorithm}
\usepackage{algorithmic}
\usepackage{booktabs}
\usepackage{nicefrac}
\usepackage[
]{todonotes}
\presetkeys{todonotes}{inline}{}
\usepackage{footmisc}

\usepackage[singlespacing]{setspace}
\newcommand{\pdiffn}[3]{\ensuremath{\frac{\partial^{#3} #1}{\partial #2^{#3}}}}
\DeclareMathOperator\erf{erf}

\newcommand{\ddiffn}[3]{\frac{d^{#3} #1}{d#2^{#3}}}     %

\DeclareMathOperator{\rhs}{rhs}

\newcommand{\vardiff}[2]{\frac{\delta{#1}}{\delta {#2}}}     %

\newdimen\CdotAxis
\newcommand*{\CdotAux}[3]{%
  {%
    \settoheight\CdotAxis{$#2\vcenter{}$}%
    \sbox0{%
      \raisebox\CdotAxis{%
        \scalebox{#1}{%
          \raisebox{-\CdotAxis}{%
            $\mathsurround=0pt #2#3$%
          }%
        }%
      }%
    }%
    \dp0=0pt %
    \sbox2{$#2\bullet$}%
    \ifdim\ht2<\ht0 %
      \ht0=\ht2 %
    \fi
    \sbox2{$\mathsurround=0pt #2#3$}%
    \hbox to \wd2{\hss\usebox{0}\hss}%
  }%
}

\def\addlegendimage{\csname pgfplots@addlegendimage\endcsname}

\DeclareSIUnit{\molpc}{mol\text{-}\%}

\title{Numerical efficiency of explicit time integrators for phase-field models}

\author[1,3]{Marco Seiz}
\ead{marco@kit.ac.jp}
\author[1,2,3]{Tomohiro Takaki}
\ead{takaki@kit.ac.jp}
\affiliation[1]{organization={Faculty of Mechanical Engineering, Kyoto Institute of Technology},
addressline={Matsugasaki, Sakyo-ku},
postcode={606-8585},
city={Kyoto},
country={Japan}}
\affiliation[2]{organization={High-Performance Simulation Research Center, Kyoto Institute of Technology},
addressline={Matsugasaki, Sakyo-ku},
postcode={606-8585},
city={Kyoto},
country={Japan}}
\affiliation[3]{Corresponding authors: marco@kit.ac.jp, takaki@kit.ac.jp}

\begin{document}

\begin{frontmatter}
\begin{abstract}
Phase-field simulations are a practical but also expensive tool to calculate microstructural evolution.
This work aims to compare explicit time integrators for a broad class of phase-field models involving coupling between the phase-field and concentration.
Particular integrators are adapted to constraints on the phase-field as well as storage scheme implications.
Reproducible benchmarks are defined with a focus on having exact sharp interface solutions, allowing for identification of dominant error terms.
Speedups of 4 to 114 over the classic forward Euler integrator are achievable while still using a fully explicit scheme without appreciable accuracy loss.
Application examples include final stage sintering with pores slowing down grain growth as they move and merge over time.
\end{abstract}
\begin{keyword}
phase-field, time integration, grain growth, final stage sintering
\end{keyword}

\end{frontmatter}

\section{Introduction}
The phase-field method is applied in many contexts to determine how microstructures evolve, with goals such as optimizing a process or explain its mechanisms.
A phase-field model typically consists of a specified free energy, which is minimized by an appropriate gradient flow, resulting in a set of partial differential equations (PDEs) for the spatiotemporal evolution.
Since these PDEs are typically not analytically solvable, numerical solutions are required and thus discretizations in space and time.
Hence the efficiency of these discretizations becomes an interesting topic to study.

For the spatial discretization a number of ``soft'' factors play a major role, such as whether an employed framework allows for adaptive mesh refinement, whether it is parallelized (efficiently), whether the order of spatial approximation can be changed easily and many more.
In contrast, temporal discretization is often done with the method of lines on top of a spatial semi-discretization.
Assuming the employed time integrator covers the relevant parts of the stability domain, the efficiency can be investigated somewhat in isolation from the specific spatial discretization and so may be transferable between different spatial discretizations.
Furthermore, if explicit schemes are employed, then no additional ``soft'' factors such as e.g. the ease of constructing appropriate preconditioners enters the question, making the decision among such schemes particularly easy.
Hence the goal of the present paper is to investigate the time integration efficiency of various explicit integrators for phase-field applications.

We note that the time evolution of different processes in the microstructure may span several orders of magnitude, making the simple forward Euler scheme often quite inefficient.
This may simply stem different timescales of the equations, but may also couple into externally given control of a field parameter.
Examples of the latter include e.g. the thermal control during solidification and sintering, where process times may be on the order of minutes to hours depending on the cooling/heating rate.
Resolving the microstructure with relevant diffusivities then often puts the stable time step on the order of microseconds or less, making full-process simulations expensive to unfeasible.
Another consideration, especially for sintering, is that the process speed depends on the grain size $G$ as a power function.
For example, if densification occurs via grain boundary diffusion the densification speed can scale as $G^{-4}$\cite{Rahaman2017}.
Hence increasing the scale of simulation from e.g. a grid spacing of 1 to 10 makes this process $10^{4}$ times slower, but a forward Euler scheme gains only a factor of $10^2$ in its stable timestep, which means simulating to equivalent density takes 100 times as long!
As the physical process occurs that much slower, unconditionally stable integrators can easily increase their time step substantially with little precision loss.

For the present paper, the focus will be on models with so-called obstacle potentials, which take on finite energy values in part of the domain, but infinite outside of it.
This causes the evolution of the phase-field to be trivial outside of a narrow-band called the interface, with the remaining part of the domain being called bulk.
Identification of the narrow-band can then be used to speed up simulations by not requiring the calculation of updates outside of it.
The time integration needs to account for the infinite potential, which is often done using a projection approach\cite{Garcke1999,Daubner2023}.

Although there is a large literature of integrators specifically crafted for the Cahn-Hilliard equation\cite{Gomez2011,Christlieb2014,Vignal2017,Xu2019,Guillen-Gonzalez2014,Cheng2021,Botchev2024} and many more, where the fourth-order nature makes the stiffness much more acute, a thorough efficiency comparison between them is more rare \cite{Christlieb2014,Vignal2017,Cheng2021,Botchev2024}; the focus is generally on assuring energy stability, i.e. that the monotonic decrease of energy predicted by the continuous problem also holds in the discrete case.
This applies even more so for the Allen-Cahn equation\cite{Vignal2017,Xu2019,Cheng2021}, with no comparative efficiency investigations being found in the literature at all for obstacle potentials with explicit integrators.
Hence this work intends to start to fill this gap by answering the following set of questions:
\begin{itemize}
 \item Does the energy stability of a time integrator play a significant role for convergence to the sharp interface limit?
 \item How coarse can integration tolerances be taken before they dominate spatial and interface width errors?
 \item How much quicker can explicit time integration become over the forward Euler integration for obstacle potentials?
\end{itemize}

To answer these questions for particular cases, benchmark geometries with analytical solutions are defined, which allow for a more appropriate performance comparison than using a reference solution.
In contrast to the benchmarks suggested by \cite{Jokisaari2017} the analytical solutions allow for exact and cheap calculation of the error with respect to the sharp interface problem, which allows the differentiation of space, time and interface width errors and is independent of the chosen free energy formulation.
Furthermore, these can be seen as extensions for benchmarking Allen-Cahn problems with optional coupling to diffusion equations without requiring 4th order Cahn-Hilliard equations to be solved.

However, these do suffer from usually being restricted to simpler geometry and parameter combinations.
Hence a final geometry, grain growth with and without mobile or immobile particles, is also considered as a more complex geometry spanning several orders of magnitudes in its parameters.
 
\section{Model}
A phase-field model using an obstacle potential will be employed within the present work, with coupling to mass diffusion realized via a Kim-Kim-Suzuki approach\cite{Kim1999}.
The independent variables are the phase-field $\phi(x, t)$ with $N$ entries and the concentration $c(x, t)$ with $K=1$ independent entries whose arguments will be dropped for brevity.
The energy functional for this reads
\begin{align}
 F &= \int f_v dV = \int a+w+f dV \label{eq:energy}\\
 a &= - \sum_{\alpha=0}^{N-1}\sum_{\alpha > \beta}^{N-1} A_{\alpha\beta} \nabla \phi_\alpha \nabla \phi_\beta\\
 w &=
 \begin{cases}
  \sum_{\alpha}^{N-1}\sum_{\alpha > \beta}^{N-1} B_{\alpha\beta}  \phi_\alpha \phi_\beta & \phi \in \Delta^{N} \\
  \infty & else \label{eq:wterm}
 \end{cases}\\
 f &= \sum_{\alpha=0}^{N-1} h_\alpha(\phi) g_\alpha(c_\alpha)
\end{align}
with the phase-field $\phi$, its individual components $\phi_\alpha$, model parameters $A, B$ which will be shortly related to physical ones, a weighting function $h_\alpha = \phi_\alpha$ and the Gibbs free energy density of each phase $g_\alpha$ in terms of the phase concentration $c_\alpha$.
The potential energy $w$ assigns infinite energy to points outside of the standard $n$-simplex
\begin{align}
 \Delta^{N} = \{ (\phi_0, \ldots, \phi_{N-1}) : \sum_{\alpha=0}^{N-1} \phi_\alpha = 1 \wedge \phi_\alpha \ge 0, \forall \alpha \in \{0, \ldots, N-1\} \}
\end{align}
in order to obtain a finite interface width.

The time evolution of $\phi$ is assumed to follow a variational derivative subject to the constraint $\sum_{\alpha=0}^{N-1} \phi_\alpha = 1$.
This constraint is built into the time evolution following \cite{Steinbach1999} by a pairwise construction
\begin{align}
 \pdiff{\phi_\alpha}{t} &= \frac{1}{N_z} \sum_{\beta=0}^{N-1} L_{\alpha\beta} n(\phi_\alpha,\nabla \phi_\alpha)n(\phi_\beta,\nabla \phi_\beta) (\vardiff{F}{\phi_\alpha} - \vardiff{F}{\phi_\beta}) \label{eq:phieq}\\
 \vardiff{F}{\phi_\alpha} &= \pdiff{f_v}{\phi_\alpha} - \nabla \cdot \pdiff{f_v}{\nabla \phi_\alpha}
\end{align}
with the pairwise mobility $L_{\alpha\beta}$ and the number of phases $N_z$ which have nonzero gradients at a point.

This construction can be thought of as a constrained gradient flow on the manifold $\sum_{\alpha=0}^{N-1} \phi_\alpha = 1$.
Note that the sum here is also only taken over phases which have nonzero gradients, i.e. whose evolution is not trivial.
To keep the sum formulation consistent with this, the function $n(\phi_\alpha,\nabla \phi_\alpha)$ is introduced which is 1 if this phase-field is nontrivial, but zero otherwise.

Reducing \cref{eq:phieq} for two phases while accounting for the sum constraint and demanding equilibrium, one obtains an analytically solvable ODE with a solution of the form
\begin{align}
\phi(x) = \label{eq:profile}
\begin{cases}
 0 & x \leq -\frac{\pi W}{2}\\
 1 & x \geq +\frac{\pi W}{2}\\
 \frac{1}{2} (1+\sin(\frac{x}{W})) & else
\end{cases}
\end{align}
with the interface parameter $W$, which relates to the interface width $L_\phi = \pi W$.
The solution together with demanding that \cref{eq:energy} is equal to the interface energy $\gamma$ when integrated over an equilibrium profile gives $A_{\alpha\beta} = \frac{4W\gamma_{\alpha\beta}}{\pi}$ and $ B_{\alpha\beta} =\frac{4\gamma_{\alpha\beta}}{\pi W}$ assuming a constant $W$.
The phase-field mobility $L$ can be related to the physical mobility by considering a pure system $f=0$ with an interface whose radius of curvature satisfies $\frac{1}{\kappa} \gg W$.
This eventually yields $L_{\alpha\beta} = \frac{\pi M_{\alpha\beta}}{4W}$ after comparing to the sharp interface solution for motion by mean curvature flow with $M_{\alpha\beta}$ being the mobility of the relevant interface.
Note that if there are driving forces besides curvature, the law of motion would change, which the formula for the phase-field mobility would need to account for in order to recover the proper kinetics.
This is generally done with an asymptotic analysis.

The concentration evolution follows a conservative Cahn-Hilliard ansatz
\begin{align}
 \pdiff{c}{t} &= -\nabla\cdot (-M_c\vardiff{F}{c})\\
 &= \nabla \cdot M_c \nabla \sum_{\alpha=0}^{N-1} h(\phi_\alpha, \phi) \frac{\partial g_\alpha}{\partial c}
\end{align}
with some mobility $M_c$.
As $c$ is the independent variable whereas the Gibbs energies are formulated in their phase concentrations $c_\alpha$, some additional relation is required.
For this the concentration is additionally represented as
\begin{align}
 c &= \sum_{\alpha=0}^{N-1} h_\alpha c_\alpha \label{eq:cadef} \\
 \mu_\alpha &= \mu_\beta = \ldots = \mu \label{eq:mueq} \\
 \mu_\alpha &= \pdiff{g_\alpha}{c_\alpha}
\end{align}
i.e. $c$ is a weighted sum of the phase concentrations and local equilibrium is assumed by equality of the chemical potential across phases. We may thus simplify
\begin{align}
 \sum_{\alpha=0}^{N-1} h(\phi_\alpha, \phi) \frac{\partial g_\alpha}{\partial c} = \mu
\end{align}
by using the chain rule and the $c$ derivative of the weighted sum representation.
Inserting back into the evolution equation yields
\begin{align}
 \pdiff{c}{t} &= \nabla \cdot M_c \nabla \mu \\
 &= \nabla \cdot \sum_{\alpha=0}^{N-1} (D_\alpha h_\alpha (\pdiffn{g_\alpha}{c_\alpha}{2})^{-1} \nabla \mu)\\
 &= \nabla \cdot \sum_{\alpha=0}^{N-1} (D_\alpha h_\alpha \nabla c_\alpha) \label{eq:ceq}
\end{align}
in which the form for $M_c = \sum_{\alpha=0}^{N-1} D_\alpha h_\alpha (\pdiffn{g_\alpha}{c_\alpha}{2})^{-1} $ is chosen such that it cancels the relevant part from $\nabla \mu = \pdiffn{g_\alpha}{c_\alpha}{2} \nabla c_\alpha$.
The last line writes out the sum and then replaces $\nabla \mu$ with its equivalent in terms of $\nabla c_\alpha$ before moving back into sum form.
One can thus interpret the concentration evolution as a weighted average of phase-diffusion problems.
The phase-wise concentrations $c_\alpha$ may be obtained from solving the system described by \cref{eq:cadef,eq:mueq}.

We employ parabolic Gibbs free energies of the form
\begin{align}
 g_\alpha = \frac{k_\alpha}{2} (c_\alpha-c_{0,\alpha})^2 \label{eq:gibbs-energy}
\end{align}
and distinguish up to two phases $\alpha, \beta$ with $\beta$ assumed to be a solid which can take on many orientations while employing the same chemical energy.
The associated grand potential is
\begin{align}
 \psi_\alpha &= g_\alpha - \mu c_\alpha = -\frac{\mu^2}{2k_\alpha} - \mu c_{0,\alpha}.
\end{align}

\section{Methods}
In this section the time integrators used to solve \cref{eq:phieq,eq:ceq} are introduced together with some additional implementation details.
The numerical framework used for the investigation supports execution on both the CPUs and GPUs, with parallelization achieved using the Message Passing Interface (MPI).
It will be detailed in another publication.
Spatial discretization is generally carried out with second order finite differences, with divergences being resolved as a sum of flux differences between faces.
Point values are laid out on a cell-centered regular Cartesian grid with spacing $\Delta x$, with each unknown in a cell representing a degree of freedom (DoF).

\subsection{Phase-field details}
In order to allow for an arbitrary number $N$ of phase-fields to be simulated efficiently, at each point in the domain only a fixed number of phase-fields $N_p$ is stored, in the spirit of \cite{Kim2006}.
This is allowable due to the finite extent of the phase-field, which takes on trivial values on a large part of the domain.
Thus each phase-field cell has $N_p$ phase-fields and integers indicating the phase these belong to, as well as an integer $N_z \le N_p$ counting phase-fields with non-zero gradients.
The number $N_p$ is adjusted to the problem at hand, with the benchmarks comparing to analytical solutions also being cheaply computable without this scheme since there $N\le 3$.
For these cases this storage scheme is unnecessary, but was still used for consistency with the final benchmark of \cref{sec:ggpores} where $N_p=12$ is employed.
The resulting runtimes for benchmarks will thus be larger than if a simple direct $N \in \{2,3\}$ scheme were to be used; the number of right-hand side evaluations however is unaffected by this.
Note that even though $N_p$ phases are formally saved, only the local phases $N_z$ are loaded and used for computation.
This is in contrast to remapping schemes such as \cite{Permann2016} which need to load and calculate the entire chosen length, even if their update is trivially zero or vanishingly small.
We furthermore note that the filter $n(\phi_\alpha,\nabla \phi_\alpha)$ in \cref{eq:phieq} is naturally adapted by this reduced storage scheme, as only nontrivial values are stored.

Due to infinite energy outside of the simplex in the energy \cref{eq:wterm} any $\phi$ that is calculated has to be within the simplex.
This is accomplished by the following trivial projection:
For any $\phi_\alpha < 0$, set its value to zero; if any $\phi_\alpha \ge 1$, set it to one and set all other $\phi_\beta$ to zero.
After this part the $\sum_{\alpha=0}^{N-1} \phi_\alpha=1$ constraint may no longer be fulfilled, which is accounted for by dividing the $\phi$ vector by its sum.
For the final geometry described in \cref{sec:ggpores} we employ a mobility-based simplex as proposed in \cite{Daubner2023}.
A more detailed investigation of the influence of the chosen simplex together with mobility contrasts is left to future work.

\subsection{Temporal discretization}
A variety of explicit temporal integrators is evaluated in the present paper, with a short overview given in the present section.
Among them is the ``standard'' forward Euler integrator as used by many explicit phase-field codes as a reference to how much faster other integrators may be.

Another integrator group is those suited for integration together with advective terms, as their stability region includes a part of the imaginary axis and they have the so-called strong stability property (SSP), the SSP schemes; due to this property, monotonic properties of the continuous time evolution can be preserved, which is critical for hyperbolic problems\cite{Ketcheson2008}.
Finally, unconditionally stable integrators built for the solution of parabolic problems are also employed, the so-called super timestepping (STS) schemes\cite{Sommeijer1998,Meyer2012,Meyer2014}.

STS schemes are chosen as they allow for more efficient time integration for semi-stiff problems and SSP schemes for potential inclusion of advective terms, which will be the focus of a later work.
SSP schemes may also be seen representative of ``normal'' Runge-Kutta methods which improve temporal accuracy but do not gain computational efficiency (no. RHS evaluations vs. stability limit) over the forward Euler integrator.

For this section we consider PDEs of the form
\begin{align}
 \pdiff{u}{t} &= f(u, \nabla u, \nabla^2 u, \ldots)\\
 \pdiff{u}{t} &\approx \rhs(u(t))
\end{align}
with a finite difference semi-discretization being carried out in transforming the continuous right-hand side $f$ into its semi-discrete part $\rhs$.

The simplest integrator considered is the forward Euler (FE, FEuler) integrator\cite{LeVeque2007} given by
\begin{align}
 u^{n+1} = u^{n} + \Delta t \rhs(u^n) + \frac{1}{2}\Delta t^2 \pdiffn{u}{t}{2} + O(\Delta t^3)
 \approx u^{n} + \Delta t \rhs(u^n)
\end{align}
with a local error of $O(\Delta t^2)$ and so global error order 1.
Its stability limit can be determined from the Dahlquist test equation and its corresponding discretization
\begin{align}
 \pdiff{u}{t} &= \lambda u \nonumber\\
 u^{n+1} &\approx (1+\Delta t \lambda)u^{n} =  R(\Delta t \lambda)u^{n}  \label{eq:stability-polynomial}
\end{align}
in which the stability polynomial $R(\Delta t \lambda)$ is introduced.
Explicit Runge-Kutta (RK) methods may always be written in the form \cref{eq:stability-polynomial}, with the only difference being the expression for $R$.
In order for a method to be absolutely stable, i.e. errors introduced earlier do not grow exponentially, the bound\cite{LeVeque2007}
\begin{align}
 |R| < 1
\end{align}
has to be satisfied, which for the above test equation (assuming $\lambda < 0$ for notational simplicity) is the case for
\begin{align}
\Delta t_e < \frac{2}{|\lambda|}.
\end{align}

We also consider two strong positivity preserving (SSP) schemes, namely SSP(n)2 and SSP(10)4, which are respectively an n-stage scheme of order 2 and a 10-stage scheme of order 4, both of which are described in \cite{Ketcheson2008}.
These schemes are efficient and accurate integrators for advective terms.
The implementation is effectively the one given in \cite{Ketcheson2008} except that an explicit additional vector is used since fulfilling the storage assumption is quite difficult on GPUs.
For example, the SSP(n)2 scheme is written in pseudocode in \cref{alg:ssp2} with $q_3$ being introduced as an additional buffer to evade the problem of reading and writing to $q_1$ concurrently.
CPU-focused codes may avoid this by computing $\rhs(q_1)$ in a layered fashion and writing back to $q_1$ as cells enter and exit the domain of dependence.
Note that for the phase-field, the simplex constraint is applied to every stage of the solution and hence the stages evolve on the appropriate sub-manifold.
The structure for the SSP(10)4 scheme is similar to the SSP(n)2 scheme and hence will not be repeated here.
\begin{algorithm}[H]
\begin{algorithmic}
\STATE $q_1 = u^{n}$
\STATE $q_2 = u^{n}$
\STATE $q_3$ // extra buffer
\FOR{$i=1$ to $s-1$}
  \STATE $q_3 = q_1 + \Delta t \frac{\rhs(q_1) }{s-1}$
  \IF{$u$ is phase-field} \STATE \texttt{simplex}($q_3$) \ENDIF
  \STATE \texttt{swap\_buffers}($q_3$, $q_1$)
\ENDFOR
\STATE $q_1 = \frac{(s-1)q_3 + q_2 + \Delta t \rhs(q_3)}{s}$
\STATE $u^{n+1} = q_1$
\end{algorithmic}
\caption{Pseudocode for the SSP(s)2 scheme using an extra buffer variable.}
\label{alg:ssp2}
\end{algorithm}

Following \cite{Ketcheson2008}, we have the stability limits
\begin{align}
 \Delta t_{SSP(n)2} &< \frac{n-1}{n}\Delta t_e\\
 \Delta t_{SSP(10)4} &< 6\Delta t_e
\end{align}
with $\Delta t_e$ being the stable FE step as derived previously.

Finally, super timestepping (STS) schemes of order 1 and 2 are also considered within the work.

These are Runge-Kutta schemes whose stability polynomial $R$ as introduced earlier is constructed in such a way to allow the stability region to grow superlinearly with the number of stages\cite{Verwer1990,Meyer2012}.
These schemes use orthogonal polynomials (e.g. Chebyshev polynomials \cite{Sommeijer1998}, Legendre polynomials \cite{Meyer2012,Meyer2014}, \ldots) as the stability polynomials $R_j$ of each stage $j$ in the resulting Runge-Kutta scheme\cite{Verwer1990}.
By relying on the recursive properties of the chosen orthogonal polynomials, these schemes can be written in a recursive form, e.g. for $s$ stages\cite{Meyer2014}:
\begin{align}
 u^{n,0} &= u^{n}\\
 u^{n,1} &= u^{n,0} + \tilde{\mu_1}\Delta t \rhs(u^{n,0})\\
 u^{n,j} &= (1-\mu_j - \nu_j) u^{n,0} + \mu_j u^{n,j-1} + \nu_j u^{n,j-2} \nonumber\\
 &+ \Delta t (\tilde{\mu_j} \rhs(u^{n,j-1}) + \tilde{\gamma_j} \rhs(u^{n,0}) )\\
 u^{n+1} &= u^{n,s} \label{eq:sts}
\end{align}
that is to say, the zeroth stage $u^{n,0}$ is the the old solution and the first stage $u^{n,1}$ is a scaled Euler step.
The remaining stages form a sequence which depends on only a limited number of prior stages to compute their stage solutions $u^{n,j}$, which means that even though hundreds of stages may be used, the used memory is limited to a small number of solution vectors, which makes it a low-storage Runge-Kutta method.
The final stage $s$ represents the solution at $t_{n+1} = t_n + \Delta t$.
The coefficients $\mu_j, \nu_j, \tilde{\mu_j}, \tilde{\gamma_j}$ depend on the chosen stability polynomial basis and the order of the scheme, with schemes of order 1 having $1-\mu_j-\nu_j=0, \tilde{\gamma_j} =0$ and hence no need to store information about the zeroth stage.
An implementation of an STS1 scheme with $s$ stages is given in \cref{alg:sts1} which uses a periodic index $j$ over a fixed number of buffers $q$.
The STS2 scheme is quite similar, except that the old solution and initial RHS are also stored and the update uses the full expression from \cref{eq:sts}.
\begin{algorithm}[H]
\begin{algorithmic}
\STATE $q$ // 3 stage buffers
\STATE $q[0] = u^{n}$
\STATE $q[1] = q[0] + \Delta t \tilde{\mu_j} \rhs(q[0])$
\FOR{$i=2$ to $s$}
    \STATE $j=i \mod 3; j_p = (i-1) \mod 3; j_{pp} = (i-2) \mod 3;$
    \STATE $q[j] = \mu_j q[j_p] + \nu_j q[j_{pp}] + \Delta t \tilde{\mu_j} \rhs(q[j_p])$
    \IF{$u$ is phase-field} \STATE \texttt{simplex}($q[j]$) \ENDIF
\ENDFOR
\STATE $u^{n+1} = q[j]$
\end{algorithmic}
\caption{Pseudocode for the STS1 scheme.}
\label{alg:sts1}
\end{algorithm}

A common feature of STS schemes is that the stable time step scales as
\begin{align}
 \Delta t = k \Delta t_e s^2 \label{eq:sts-limit}
\end{align}
i.e. quadratically in the number of stages $s$ relative to the stable FE step (largest eigenvalue).
Thus these schemes can be said to be unconditionally stable, as one might simply choose a large number of stages to achieve the desired stability limit.
Choosing a larger step also makes the scheme more efficient relative to the FE scheme as increasingly fewer RHS evaluations will be required; this is in contrast to the SSP(n)2 scheme, whose stability domain only grows linearly.
Within this work we only consider the Legendre basis introduced by \cite{Meyer2012,Meyer2014} for simplicity, with schemes thus becoming Runge-Kutta-Legendre (RKL) schemes; the general trends should follow for other bases, with a possible exception of a Chebyshev basis when using spatially variable coefficients, as these may introduce step artifacts in the solution\cite{Meyer2014}.
In order to estimate the stable Euler time step $\Delta t_e$ necessary to compute the stable number of stages via \cref{eq:sts-limit}, we employ Gershgorin's circle theorem on \cref{eq:phieq,eq:ceq} separately for a two-phase system.
This results in estimates $\lambda_\phi, \lambda_c$ which are detailed in the supplementary material.
The larger of the two resulting eigenvalues is then used to estimate a stable Euler step $t_e = \frac{2}{\lambda}$, which is multiplied by a safety factor of \num{0.9} before being used to estimate the stage number, which is rounded to the next largest odd number as even stage numbers fail at damping short wavelengths\cite{Meyer2012}.
If a problem only involves the phase-field, as in \cref{sec:graingrain}, only the phase-field eigenvalue estimate is employed without terms relating to the chemical driving force.
For simplicity no nonlinear power method was used for estimating the eigenvalues; for more complicated couplings in which Gershgorin's circle theorem yields very crude bounds or involves substantial linearization it may also speed up the computation by reducing the necessary stage count.

The simplex constraint of the phase-field suggests one deviation from the STS2 scheme:
If in the first stage the new phase-field values exceed the simplex constraint, adding the same RHS again later will most likely continue to trigger the simplex, introducing additional error.
Thus, if this case is detected, then $\rhs(u^{n,0})$ is re-computed as
\begin{align}
 \rhs(u^{n,0}) &= \frac{S(u^{n,1}) - u^{n,0}}{\tilde{\mu_1}\Delta t}
\end{align}
with $S$ being the simplex projection operator.
The effective initial RHS is one s.t. the simplex constraint is not violated.
In computational experiments this significantly improved the error behaviour.

\subsection{Temporal adaptivity}
For conditionally stable schemes, often a time step close to the stability limit is chosen for reasons of computational efficiency, with the hope being that the temporal error introduced thus is much smaller than the spatial error.
However, for unconditionally stable schemes some kind of adaptivity is necessary as to decide what step size to take.
For this an estimate of the local error introduced by the time stepping process is useful.
The SSP methods may be extended with an embedded method, i.e. a method requiring function evaluations at the same or a subset of stages, and which calculates a lower-order solution of the problem.
The difference between the solutions may then be taken as a measure of the error, i.e. $e^{n}_i = u^{n}_i - z^{n}_i$ with $z$ being the embedded solution and $u$ being the normal solution.
Hence embedded methods can estimate the error at little extra computational cost but with one extra vector being stored.
The embedded methods for the SSP methods use the embedded solution coefficients suggested by \cite{Fekete2022}.

For the STS methods, the error may be estimated following \cite{Sommeijer1998} by a calculation involving the old and new solution.
This requires an extra RHS evaluation before a time step is accepted\footnote{This evaluation might be re-used for the first stage of the next step if it's accepted; for simplicity this wasn't implemented and the additional RHS evaluations were counted.} but without extra storage requirements.
The error formula is then
\begin{align}
 e^{n}_i &= \frac{1}{15}(12 (u^{n+1}_{i}-u^{n}_i) + 6 \Delta t (\rhs(u^{n}_i) + \rhs(u^{n+1}_{i})))
\end{align}
which effectively calculates a scaled solution estimate for step $n+1$, which for the phase-field needs to be subject to simplex constraints.
In order to account for this, the phase-field error components are locally filtered:
Any component of a trial solution $\tilde{u}(t+\Delta t) = u^{n} + \frac{\Delta t}{2} [\rhs(u^{n}) + \rhs(u^{n+1})]$ which lies outside of $[0,1]$ is assumed to have zero error.
Once a phase-field value leaves this interval, it is reasonable to assume that it will not immediately enter it again, as the phase-field describes the advance of a free boundary.

In either case, once a local error estimate $e^{n}_i$ is available, the error of this step is calculated as\cite{Sommeijer1998}
\begin{align}
 e^{n} &= \left[\frac{1}{N} \sum_{i=0}^{N} |\frac{e^{n}_i}{r \max(u^{n}_{i}, u^{n+1}_{i}) + a}|^2\right]^{\frac{1}{2}} \label{eq:error_norm}
\end{align}
by computing a weighted $l^2$ norm of the error over all DoFs $i$.
The weight is calculated with a relative tolerance $r$ and an absolute tolerance $a$ to define the level of acceptable error.
These tolerances are defined per independent field, i.e. $r_\phi, a_\phi$ for the phase-field and $r_c, a_c$ for the concentration.
One would usually choose the absolute tolerance at the level of error or noise one expects of a field, and so it may be seen as a threshold at which a value might as well be zero.
For the concentration $c$ this lets us choose reasonable tolerances since we know the expected value range via thermodynamics and its smallest reasonable value doesn't approach zero.
For the phase-field however the entire interval $[0,1]$ represents possible values, with one cell often containing values close to both $1$ and $0$.
Since the phase-field can get arbitrarily close to $0$ and in fact can become $0$, the choice of $a_\phi$ is not obvious.
Hence it will be varied within the work and its influence mentioned where appropriate.
While assigning different tolerances to the different entries of the phase-field is possible, each individual phase-field can range over $[0,1]$; only the sum is known to be 1.
Future work may consider altering how the error is calculated for the phase-field.

If $e^{n} < 1$, the step is accepted, otherwise rejected.
In either case, a new time step is calculated following the PID-controller of \cite{Soderlind2006} together with a bias $b \in [0.1, 0.98]$ and the order $P$ of the method
\begin{align}
 E_{n} &= e^{n}/b\\
  F &= (E^{-k_1}_{n} E^{-k_2}_{n-1} E^{-k_3}_{n-2})^{1/(P+1)} (\frac{\Delta t_n}{\Delta t_{n-1}})^{k_4}  (\frac{\Delta t_{n-1}}{\Delta t_{n-2}})^{k_5} \\
 \Delta t_{n+1} &= \Delta t_{n} (1 + 5 \arctan( \frac{F-1}{5}))
\end{align}
with $k_1 = 1.25, k_2 = 0.5, k_3=-0.6, k_4=0.25, k_5=0$; note that the step index $n$ is now a subscript as not to be confused with the exponentiation.
The bias $b$ starts at $0.9$ and makes the step estimate slightly more conservative; it is multiplied by itself on step rejection to make successive step rejections more unlikely, up to a value of $0.1$.
On step acceptance, it is divided by $\sqrt{b}$ but limited to $0.98$.
The time step factor $F$ is then smoothed to limit the overall step size changes.
Terms for which there is no data yet are dropped from the above expression.
We observed that the full controller may suggest increasing the time step upon rejection; in this case only $e^{-k_1/(P+1)}_{n}$ is used for estimating the new timestep, which is certain to decrease the time step.

The compact support of the phase-field and its storage scheme impacts the definition of $N$ within \cref{eq:error_norm}.
For a vectorial unknown $X$ with $k$ components, $N = Vk$ trivially with $V$ being the number of points in the domain.
Using this definition would massively underestimate the error, as most entries of the phase-field will attain trivial values and hence have zero error.
Hence we only count nontrivial entries, which are defined as having nonzero $e_i^n$.
Furthermore, due to the compact support there may be entire regions for which there are zero nontrivial DoFs.
In order to stay somewhat comparable to well-types of potentials, a zero nontrivial phase-field cell is treated as 2 DoF, since a direct two-phase implementation would simply have two DoFs per cell.
In these cases one would have a nonzero but vanishingly small error in these regions, which should be somewhat comparable to the present approach.
A more detailed investigation, e.g. only accounting for errors for the phase-field in the interface, is left to future work.

We note that reductions, such as required by the error estimator, are nondeterministic on GPUs\cite{Shanmugavelu2024} and so the runtimes and RHS counts can slightly fluctuate when using an adaptive time step on GPUs.
While this also applies to MPI parallel reductions, the number of non-associative operations is often far smaller and so can usually be ignored.

\subsection{Error sources}
\label{sec:errsource}
There are many kinds of errors when integrating a PDE numerically; this section shall serve as a small overview of the errors considered within the work.
First there are discretization errors, in both the space and time dimension, induced by the discretization schemes.
Their order can generally be determined analytically, with a concrete value requiring knowledge of the solution or computational experiments as they relate to the derivatives of the solution.
We will generally consider both.
Note that spatial and temporal error can compensate each other and thus conditions which should have larger temporal error can in fact end up with smaller total error.
Next, when working with floating point arithmetic there is necessarily some amount of error involved in mapping the infinite real numbers into a finite set.
We generally do not consider these to be significant, as the quantities we compute are scaled to be $O(1)$ such that the floating point error should be far below the discretization error.
Finally, an error specific to phase-field methods is that originating from the finite interface width, which manifests in a deviation from the target free-boundary problem.
When only comparing to reference solutions at refined space/time grids, this error does not play a role, but is critical for matching the physics of a problem.
Generally, lower interface widths produce less error, but the spatial discretization interacts with the interface width and hence some care is required when performing grid convergence studies and comparing against free-boundary problems.

A central finite difference discretization for the second derivative of a function $u(x)$ may be obtained by adding Taylor expansions for $u(x_0+\Delta x)$ and $u(x_0-\Delta x)$ and solving for the second derivative at $x_0$.
This yields, noting that all odd derivatives cancel,
\begin{align}
 \ddiffn{u}{x}{2} &\approx \frac{u(x_0+\Delta x) - 2u(x_0) + u(x_0-\Delta x)}{\Delta x^2} - 2\sum_{n=2}^{\infty} \frac{\Delta x^{2n-2}}{(2n)!} \ddiffn{u}{x}{2n}
\end{align}
i.e. the standard second order finite difference, with the infinite sum at the end representing the error.
Combine this with $\ddiffn{\phi}{x}{n} \propto \frac{1}{W^n}$ for the equilibrium profile \cref{eq:profile} away from the bulk via simple differentiation, shift indices, and we can see that the spatial discretization error scales as
\begin{align}
 e_{\Delta x} \propto \sum_{n=1}^{\infty} q_n \frac{\Delta x^{2n}}{W^{2+2n}}
\end{align}
if we are sufficiently far from the bulk to ignore the case statement, with $q_n$ absorbing any remaining constants and factors depending on $n$ solely.
Thus, in order for the error to vanish as $\Delta x, W \to 0$, the terms $\frac{\Delta x^{2n}}{W^{2+2n}}$ need to go to zero.
For constant $W$, one obtains the usual error order if a reference solution using a sufficiently small $\Delta x$ is employed.
However, this does not necessarily recover the sharp interface limit $W\to0$ of the free-boundary problem.
Thus, when doing grid convergence studies, $W$ is chosen as $W=3\Delta x^{m}$ with $m = 0.4 < \frac{1}{2}$ for which the above error goes to zero for any $n$.
This scaling essentially means that the interface becomes increasingly well-resolved as the grid spacing $\Delta x$ is reduced.
The prefactor is chosen such that at the coarsest resolution employed we have about 9 cells in the interface.
Finally, note that the error w.r.t. the sharp interface solution does not necessarily follow a simple power law, since it is the combined error of discretization and from the diffuse interface.
 
\section{Common conditions \& parametrization}
In this section the common parameters and conditions will be explained, with the following section presenting single benchmarks and their solution together.
The benchmarks are generally chosen such that an analytical solution for the equivalent free-boundary problem exists.
This allows us to determine when spatial or temporal error dominate as well as what kind of tolerances can be appropriate.
If not mentioned otherwise, all boundary conditions are periodic.

In most cases, the initial conditions of the phase-field are based on calculating the profile function \cref{eq:profile} with a suitably defined distance field $d$ as its argument.
Thus only the distance field $d$ will be given for these cases.
However, for some configurations the geometric description based solely on distance fields breaks down as there are multiple overlapping distance fields, whose summed phase-field value do not sum to 1; this is the case for \cref{sec:triplejunction,sec:ggpores}.
In this case we consistently proceed as follows:
We pick one phase $\alpha$ to exactly follow the distance field $d_\alpha$ and thus we can fix its value $\phi_\alpha = \phi(d_\alpha)$.
Based on the distance fields of the remaining grains, we generate trial $\phi$ values whose ratios $q_{ij} = \frac{\phi_i}{\phi_j}$ we would like to preserve together with the sum constraint $1-\phi_\alpha = \sum_{\beta=0, \beta \neq \alpha}^{N-1}\phi_\beta$.
This generates at first a nonlinear system for the new values $\phi^{'}$ which satisfy these constraints.
By multiplying all equations but the sum one by the unknown divisors we obtain a linear system, e.g. for $\alpha=0$:

\begin{align}
 \begin{bmatrix}
 1      &  -q_{12} &       0 & \cdots & 0\\
 0      &      1 & -q_{23} & \cdots & 0 \\
 \vdots & \vdots &  \vdots & \vdots & \vdots\\
 1      & 1      & 1       & \cdots & 1
\end{bmatrix}
\begin{bmatrix}
 \phi^{'}_1\\
 \phi^{'}_2\\
 \vdots\\
 \phi^{'}_{N-1}\\
\end{bmatrix}
=
\begin{bmatrix}
0\\0\\ \vdots\\1-\phi_\alpha
\end{bmatrix}.
\label{eq:distfix}
\end{align}

If a simulation includes \cref{eq:ceq}, then the initial concentration field is given by the local equilibrium \cref{eq:cadef} unless mentioned otherwise.
We note though that resolving a discretized profile with different interface widths tends to produce slightly different areas after integration, as there is some error associated with the interface width (see the supplementary material for details).
In order to counteract this, the difference of the sharp-interface concentration $C = V_e c_e + V_E c_E$ (with the sharp interface volumes $V_e$ and $V_E$ of the embedding/embedded phases) to the concentration in the domain $\int_V c dV$ is used to shift all concentration values in the domain by their difference divided by the volume of the domain.
This is only done for the equilibrium benchmarks \cref{sec:phaseembed,sec:triplejunction} to ensure approximately the same equilibrium sizes are recovered under grid refinement.

As equilibrium states should be independent of the choice of kinetic parameters, ignoring discretization error, unconditionally stable schemes can speed up convergence somewhat:
To a first approximation, the time to equilibrium depends on how fast a homogeneous chemical potential is reached, which scales with the driving forces and diffusivity.
The phase-field mobility only enters the problem insofar as moving the interface, following the diffusion field.
Thus increasing the diffusivity can reduce the time required to reach equilibrium.
Note that if all kinetic parameters are increased by the same factor, nothing is gained as this simply rescales the time for the system.

The values used for equilibrium determination in the benchmarks \cref{sec:phaseembed,sec:triplejunction} may oscillate around a steady-state value.
In order to account for this, an oscillation resistant-algorithm detailed in the supplementary material is employed.

The standard set of parameters employed is shown in \cref{tab:params}, which will be used unless mentioned otherwise; these are nondimensional values and do not represent any specific material.
The units they would have in dimensional form are listed together with them.
As the concentration is measured in mole fraction, it is naturally dimensionless.

For equilibrium simulations, $D=100$ is used to speed up convergence.
When performing grid convergence studies $W=3\Delta x^{0.4}$ is used as described in \cref{sec:errsource} to ensure we converge to the sharp interface limit.

\begin{table}[h]
\centering
\caption{employed parameters}
\label{tab:params}
\begin{tabular}{lll}
parameter   & value  & units\\
\hline\\
$M_\phi$ & 1 & $\si{m^4/(Js)}$\\
$\gamma$ & 1 & $\si{J/m^2}$ \\
$W$ & 2.5 & $\si{m}$ \\
$\Delta x$ & 1 & $\si{m}$\\
$D$ & 1 & $\si{m^2/s}$\\
$k_\alpha = k_\beta$ & 500 & $\si{J/m^3}$ \\
$c_{0,\alpha}$ & 0.02 & - \\
$c_{0,\beta}$ & 0.98 & - \\
\end{tabular}
\end{table}

A short summary of the integrators used in this work together with how they are run is given here; a tabular overview is given in \cref{tab:ints}.
The FEuler and SSP(5)2 schemes are run at fixed time steps, limited by their respectively stability limits.
The SSP(10)4 scheme at tolerances of \num{1e-14} for all fields is used as a reference solution generator for testing convergence order for some problems.
STS schemes of order 1 and 2 are run with both fixed time steps and adaptively.
The time steps are based on being a multiple of the stable Euler time step up to a predetermined number.
When integrating adaptively, we generally use a relative tolerance of $r=\num{1e-4}$.
The absolute tolerance for the concentration is $a_c=\num{1e-4}$ unless mentioned otherwise; this represents about 1\% of the smallest expected concentration value.
As mentioned above, the phase-field tolerance will be varied between $\num{1e-2}$ and $\num{1e-4}$ and its influence mentioned where appropriate.
When performing grid convergence studies, only the STS2 integrator with $a_\phi \in \{\num{1e-2}, \num{1e-4}\}$ will be used.

\begin{table}[h]
\centering
\caption{Employed integrator configurations}
\label{tab:ints}
 \begin{tabular}{ll}
 integrator   & configurations \\
\hline\\
FEuler & $\Delta t \in \{0.5, 1\} \cdot \Delta t_e$\\
SSP(5)2 & $\Delta t \in \{0.5, 1, 2, 4\} \cdot \Delta t_e$\\
SSP(10)4 & $a_{\phi,c} = \num{1e-14}, r_{\phi,c} = \num{1e-14}$ \\
STS[12] & \begin{tabular}{@{}l@{}}
 $\Delta t \in \{1, 10, 100, 200\} \cdot \Delta t_e$\\
 $a_{\phi} \in \{\num{1e-2}, \num{1e-3}, \num{1e-4}\}$
          \end{tabular} \\

\end{tabular}
\end{table}
 
\section{Results \& discussion}
In this section the geometry of the benchmarks together with their solution is presented.
The first subsection will deal with equilibrium problems to address questions of energy stability, with the second one dealing with kinetic benchmarks to assess the computational efficiency of the schemes.
The final subsection uses a geometry without analytical solution to show how the small-scale benchmark results translate into more practical applications.
Performance is reported here in terms of RHS evaluations, as this measure is independent of the machine employed and hence eases future comparison with other time integrators.
Qualitative comparison to implicit integrators can also be done with this as long as the resulting systems are solved with iterative methods as each iteration requires RHS evaluations; a more quantitative comparison needs to account for e.g. additional matrix-vector operations required by the solution algorithm as well as potential mass matrix inversions.
The supplementary material includes data on the runtimes as well for the interested reader, though the conclusions do not change if performance is measured in runtime instead of RHS evaluations.

Finally, we generally output at defined target times exactly, with these being available for each configuration in the supplementary data.
For fixed-step integrators this only manifests in an extra RHS evaluation for each output time not being an integer multiple of the time step width.
For adaptive integrators, this generally reduces the step taken prior to reaching the output significantly and hence becomes less efficient than if a simple time crossing output strategy were adopted.
Hence the performance results are slightly biased against adaptive integrators.

\subsection{Equilibrium verification}

An important part of the phase-field theory is that the derived PDEs represent gradient flows of an energy functional.
This ensures that the evolution minimizes that functional.
Schemes ensuring this are generally referred to as energy stable\cite{Gomez2017}, with many classical schemes not being unconditionally so.
In lieu of proving (conditional) energy stability for schemes employed herein, two cases with analytical equilibrium state are considered.
Even if the schemes should temporarily have increasing energy (as shown), it is deemed acceptable if the analytical equilibrium state is the convergence point of a grid refinement study and that at a practical resolution, the scheme approximates the equilibrium state closely as well.
Note that all schemes employed here contain at least one (scaled) forward Euler step and hence would at best be conditionally energy stable.
We further note that even energy stable schemes might not lead to the correct kinetics\cite{Xu2019}, so they should also be benchmarked before use.

\subsubsection{Phase embedding}
\label{sec:phaseembed}
The first equilibrium case is that of one phase being embedded in another, with curvature providing a driving force for shrinkage which is counteracted by chemical driving forces.
Here we solve both \cref{eq:phieq} and \cref{eq:ceq}.
Both phases are initialized with their flat equilibrium concentration values and thus there is a small driving force for diffusion in order to reach equilibrium in which there is a pressure difference, the Laplace pressure,
\begin{align}
  \Delta P = \gamma\kappa
\end{align}
between both phases, with the curvature $\kappa$\cite{Balluffi2005}.
This pressure difference is achieved by a shift in chemical potential via mass transfer.
In order to accommodate this mass transfer, the volume of the individual phases changes slightly.
The curvature is determined from a shape assumption, e.g. for a circular embedding we have $\kappa = \frac{1}{r}$ with $r$ being the radius of the embedded phase.
The pressure difference is equivalent to the difference in grand potentials $\Delta \psi = \Delta P$ of the bulk phases in equilibrium, as this is the driving force resisting capillary pressure.
The bulk grand potentials are defined as the values the respective grand potentials take on at the average chemical potential of the domain.
Hence the error is defined as
\begin{align}
 e &= |\gamma\kappa - \Delta \psi|.
\end{align}

The distance field for the inner phase is
\begin{align}
 d = r - ||X - C||_2 \label{eq:centersphere}
\end{align}
with the position vector $X$, the center of the system $C = (\frac{L_1}{2}, \ldots)$ and the initial radius of the phase $r$.
The outer phase follows as the complement.

We first perform a grid convergence study with this geometry, resolving a domain of $128^2$ with an initial seed of size $r=32$.
Equilibrium for this case is defined as a rate of change for the Laplace pressure $|\frac{(\Delta P(t) - \Delta P(t_b))} { t - t_b }|$ being smaller than \num{1e-14} for two times $t > t_b$ at least two diffusion times $128^2/D$ apart.

\begin{figure}[h]
 \centering
  \begin{subfigure}[t]{0.4\textwidth}
    \includegraphics[width=\textwidth]{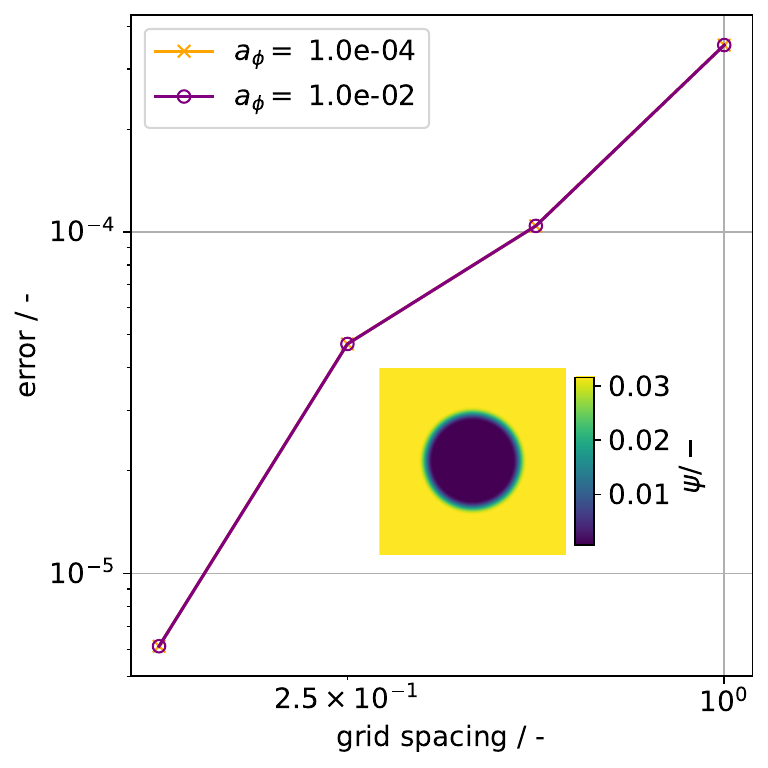}
    \caption{errors in the Laplace pressure}
    \label{fig:phasephase-convergence}
  \end{subfigure}
  \begin{subfigure}[t]{0.4\textwidth}
    \includegraphics[width=\textwidth]{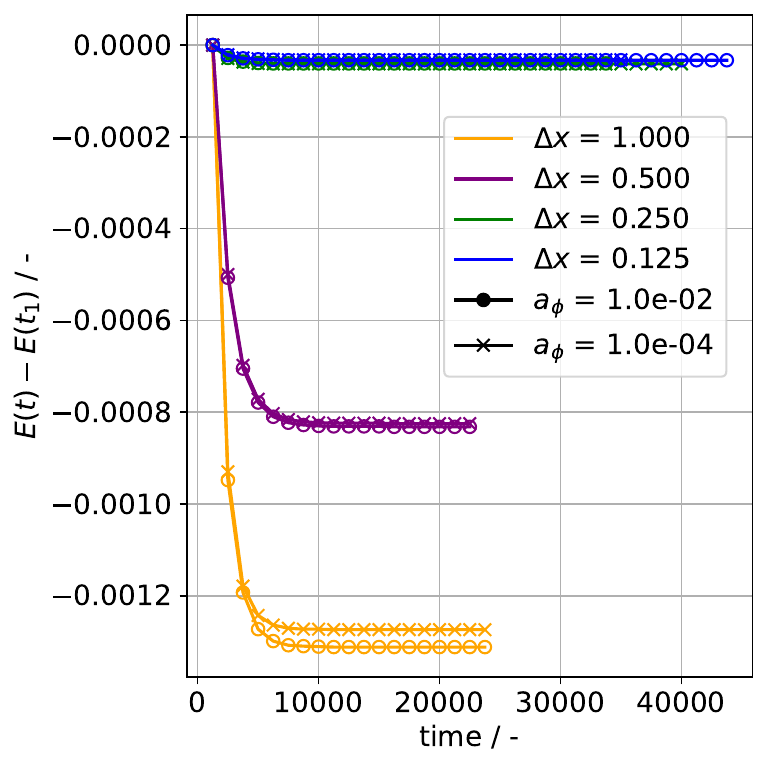}
    \caption{energy profile over time}
    \label{fig:phasephase-energy}
  \end{subfigure}
 \caption{Results for the embedded phase case, showing non-uniform convergence to the sharp interface limit. The equilibrium grand potential $\psi$ is shown exemplary as an inset; it clearly divides the two bulk phases with homogeneous values inside. The energy is observed to decrease monotonically.}
 \label{fig:phasephase}
\end{figure}

The results of this study in terms of the error and energy evolution are shown in \cref{fig:phasephase}.
The inset shows the spatial distribution of the grand potential $\psi$ in equilibrium, showing a clear delineation between the two phases.
We observe convergence, though not at a uniform rate.
At the coarsest resolution, a relative error of about $1\%$ ($\Delta P_{eq} \approx \num{0.031}$) is obtained, suggesting the equilibrium can also be well-approximated with coarse resolutions as often employed in phase-field studies.
Note that the results for the two different absolute tolerances $a_\phi$ effectively overlap, suggesting the time integration error is negligible.
For this case a monotonic decrease in energy is observed, with the energy change after the first output frame being very small; hence in \cref{fig:phasephase-energy} only the change in energy relative to this frame is shown.
0 to 2\% of steps were rejected, though given that about 50 steps were calculated, this is likely just the initial adjustment.

\subsubsection{Double triple junction}
\label{sec:triplejunction}
\begin{figure}[h]
 \centering
  \begin{subfigure}[t]{0.48\textwidth}
    \includegraphics[width=\textwidth]{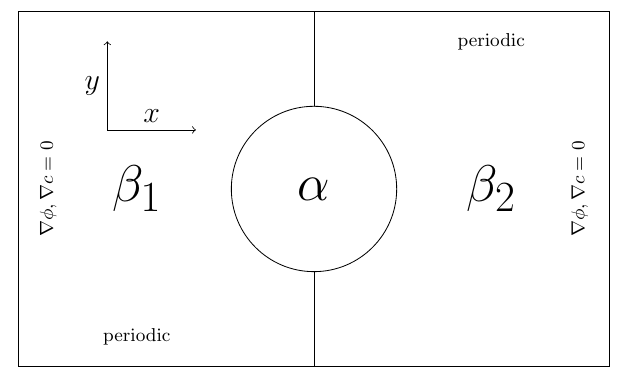}
    \caption{initial \& boundary conditions}
    \label{fig:tj-icbc}
  \end{subfigure}
  \begin{subfigure}[t]{0.49\textwidth}
    \includegraphics[width=\textwidth]{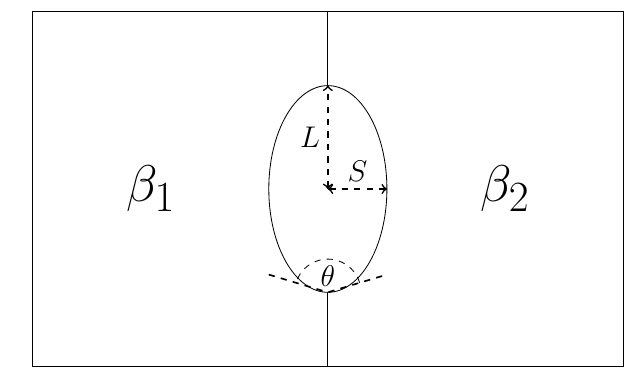}
    \caption{equilibrium geometry}
    \label{fig:tj-eq}
  \end{subfigure}
 \caption{Geometry for the double triple junction.}
 \label{fig:tj}
\end{figure}

The second case is a triple junction geometry, now also probing variable parameters and multiple phases interacting.
Again, both \cref{eq:phieq} and \cref{eq:ceq} are solved in this case.
One phase $\alpha$ is placed on the grain boundary between two grains of the phase $\beta$, as sketched in \cref{fig:tj} together with the employed boundary conditions and geometric quantities, forming two triple junctions.
The junctions will seek to equilibrate the surface tensions and hence evolve the structure\cite{Balluffi2005}.
In contrast to \cite{Daubner2023}, we want to test the equilibrium together with a chemical driving force, and so do not require Dirichlet boundary conditions to prevent curvature minimization from eliminating all interfaces.
We note that this choice generally leads to much longer simulation times, as the equilibrium is now linked to the long-range diffusion field.
The evolution should be such that the phase $\alpha$ obtains an equilibrium shape defined by the dihedral angle
\begin{align}
 \theta_{eq} = 2 \arccos(\frac{\gamma_{\beta\beta}}{2\gamma_{\alpha\beta}})
\end{align}
with its momentary form being computable with the geometric parameters $L$, $S$ following \cite{Bradley1977} by assuming the shape is a vesica piscis:
\begin{align}
 \theta = 4\arctan(\frac{S}{L}).
\end{align}
These parameters are estimated based on the extremal extents of the $0.5$ level set of phase $\alpha$.
The error is then simply defined as
 \begin{align}
  e &= |\theta_{eq}- \theta|.
 \end{align}

The phase $\alpha$'s distance field is given by \cref{eq:centersphere}, with those of the two $\beta$ grains being $d_{\beta_1} = C_x -x$ and its complement respectively.
In this case the sum of the phase-fields do not sum to 1 and hence the system \cref{eq:distfix} is solved while fixing $\phi_\alpha$ and determining new values for both $\phi_\beta$.

Here we resolve a domain of size $192 \times 96$ at the coarsest scale with $\Delta x = 1$, with the long side being along the grain boundary to facilitate simulation of smaller dihedral angles at the same aspect ratio.
The initial radius $r$ of the circular phase on the grain boundary is $32$.
For simplicity the convergence study is limited to a single angle $\gamma_{\beta\beta}=1, \gamma_{\alpha\beta}=2 \rightarrow \theta \approx \SI{151}{\degree}$, as results should be transferable to other angles except when approaching $\theta \to 0$ as complete dewetting is not captured by the model \cite{Daubner2023}.
In addition, the coarsest resolution is also run with the FEuler integrator with a time step 0.9 times its stable time step.
The remaining angles in the range \SIrange{95}{180}{\degree} will use the practical resolution $\{\Delta x = 1, W = 2.5\}$ with $\gamma_{\alpha\beta} \in \{0.75, 1, 2, 6\}$.

Equilibrium is defined as a simple approximation to the measurement described in \cref{sec:triplejunction}.
The positions of the 0.5 lines are determined on the horizontal and vertical lines passing through the center of the domain, using only a linear interpolation in the respective coordinate directions, with their respective differences giving approximations to $L, S$.
This allows an estimate of the current dihedral angle $\theta$, with equilibrium defined as a rate of change in angle by less than $\SI{1e-11}{\radian}$ over at least two diffusion times.
The reported angle values are based on a more accurate marching squares reconstruction of the 0.5 contour lines of the phase $\alpha$, which is only done after the simulation is completed.

\begin{figure}[h]
 \centering
  \begin{subfigure}[]{0.4\textwidth}
    \includegraphics[width=\textwidth]{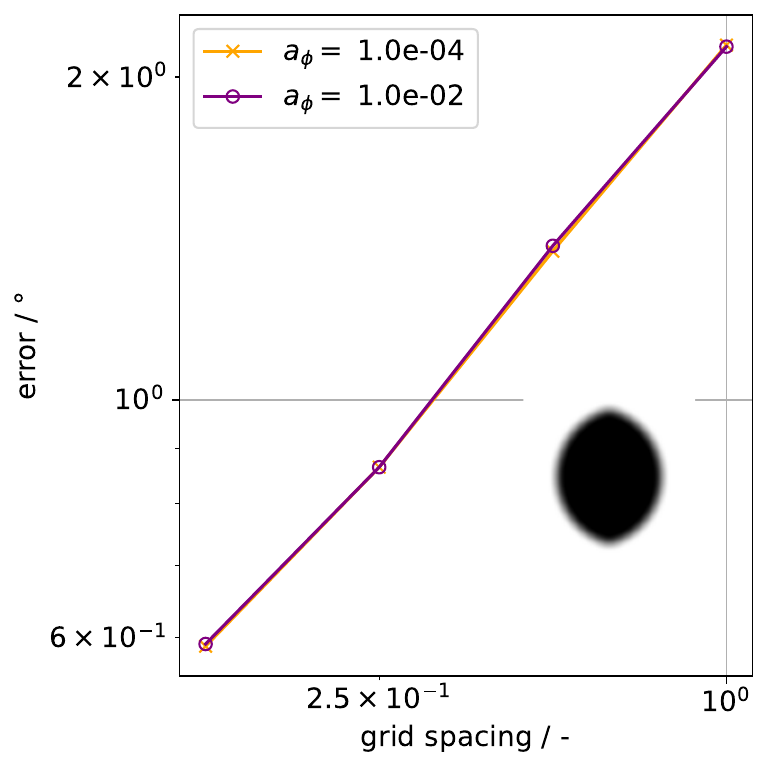}
    \caption{errors in the dihedral angle}
    \label{fig:tjangle-convergence}
  \end{subfigure}
  \begin{subfigure}[]{0.4\textwidth}
    \includegraphics[width=\textwidth]{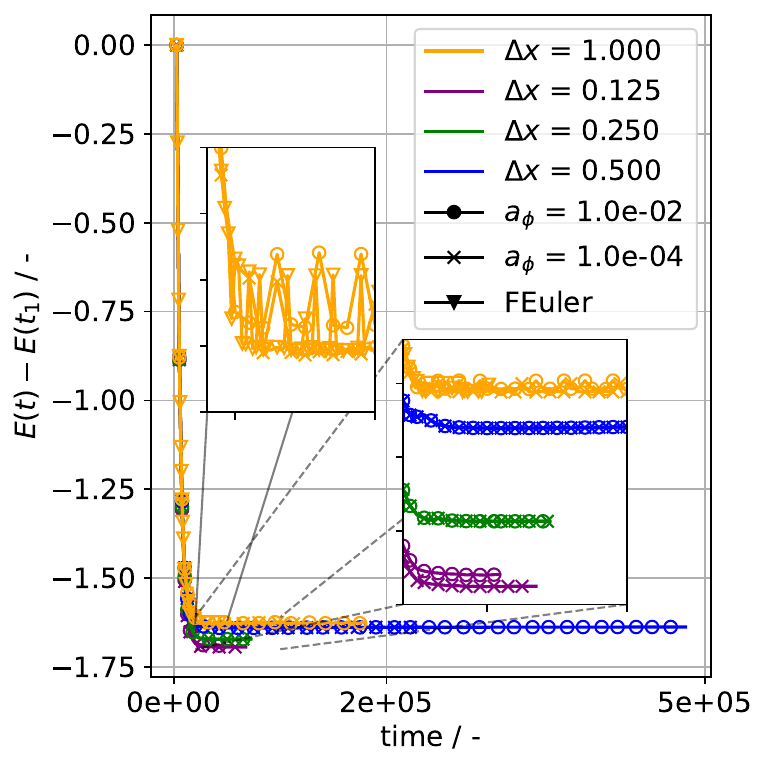}
    \caption{energy profile over time}
    \label{fig:tjangle-energy}
  \end{subfigure}
  \begin{subfigure}[]{0.8\textwidth}
    \includegraphics[width=\textwidth]{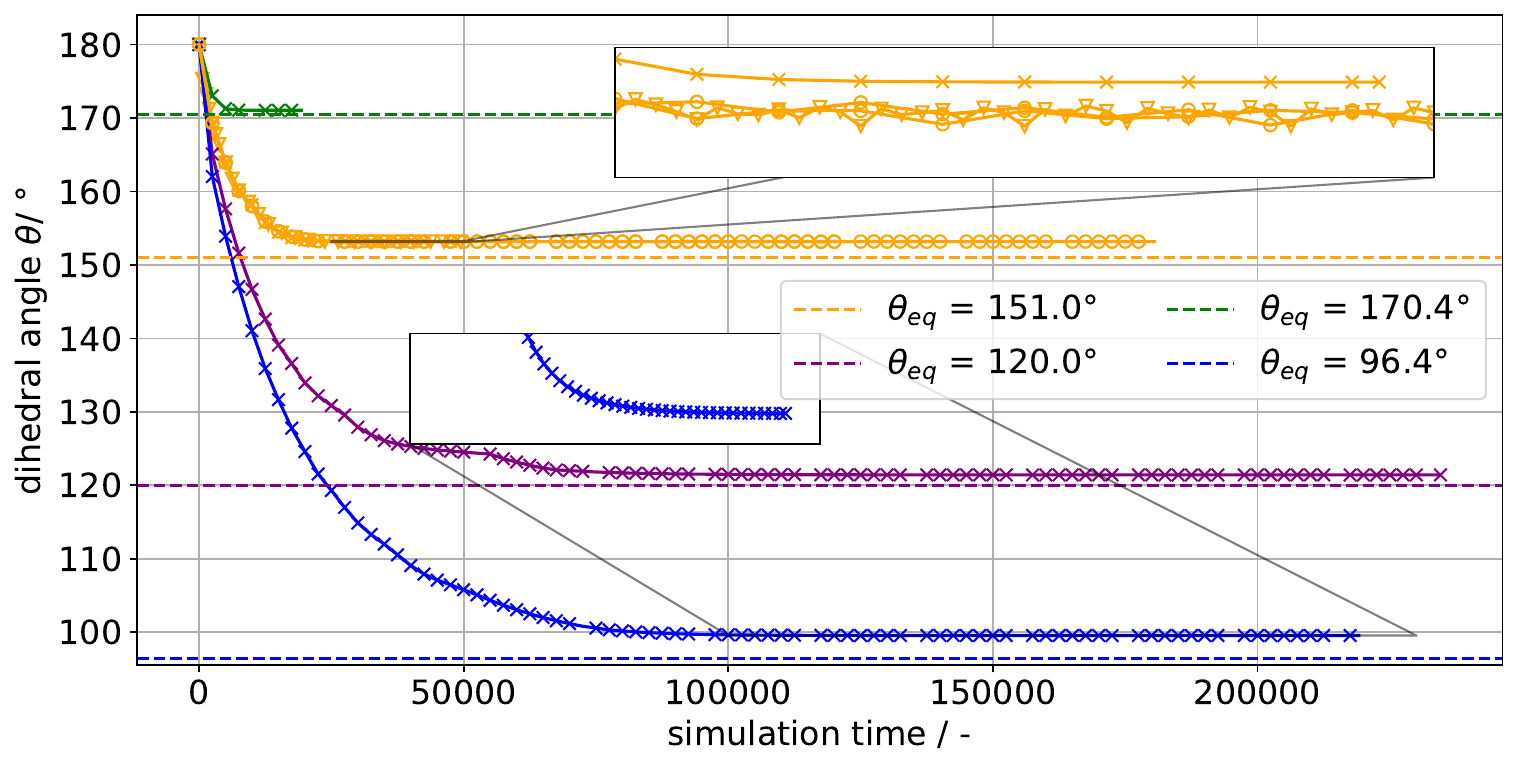}
    \caption{time evolution of the angle}
    \label{fig:tjangle-time}
  \end{subfigure}
 \caption{Results for the double triple junction case, showing convergence to the sharp interface but also non-monotonic decrease of energy; the equilibrium shape in the inset is that of the predicted vesica piscis shape. At the practical interface resolution, the analytical angle is still approximated well within the entire range of tested interface energy ratios. The circular marks in \cref{fig:tjangle-time} indicate STS2 results for $W=3$ and the triangles for the FEuler integrator.}
 \label{fig:tjangle}
\end{figure}

The results of this study are shown in \cref{fig:tjangle}.
Again, the grid convergence study in \cref{fig:tjangle-convergence} shows convergence to the sharp interface with little influence from $a_\phi$.
At the coarsest resolution, a relative error of $1.4\%$ is obtained.
In this case, the energy does not decrease monotonically beyond the first output frame as shown in the inset of \cref{fig:tjangle-energy}, though there is a long-term decrease of energy.
As is visible from the graph, on the scale of the energy reduction it is hardly noticeable without zooming closely to the data.
Note that the FEuler integrator at the coarsest resolution also exhibits a non-monotonic behaviour as shown in the top inset.
We note though that energy instability was only apparent once the simulation was close to equilibrium; benchmarks not getting close to equilibrium are unlikely to probe energy stability properties of the integrator.
If desired, one way of avoiding energy instability would be adding a second criterion for step rejection dependent on whether the energy decreased or not.
If the energy increased, integration tolerances would be tightened and a new time step computed based on this until the energy no longer increases.
This is based on the thought that error in energy is bounded by the change in the solution per step, with tightening tolerances reducing the change and hence reducing the energy error.

Results for all investigated angles are shown in \cref{fig:tjangle-time}, with errors also on the order $O(\SI{1}{\degree})$, with a tendency for the error to get larger at lower dihedral angles, as also observed by \cite{Daubner2023}.
The insets show zooms on the evolution data, revealing that not necessarily a steady state is reached, but rather that the simulation may oscillate around one.
The oscillations here are obtained for $W=3$ (circles), i.e. the coarsest part of the grid convergence study, with either employed $a_\phi$, but apparently not for $W=2.5$, though both exhibited energy instability.
Since energy instability does not strictly predict the oscillatory behaviour in the angle, it is hard to say whether these are related, but in any case the time integration error is small compared to the spatial error as evidenced by the close overlap for $W=3$.
Based on this, it seems that entering the energy instability regime has no practical influence on simulations approximating the sharp interface equivalents.

\FloatBarrier

\subsection{Kinetic verification}
Given that energy stability seems to play no practical role for equilibrium, we proceed with the kinetic benchmarks to determine the computational efficiency of the integrators.
\subsubsection{Single grain growth}
\label{sec:graingrain}
A single grain embedded within another grain will shrink due to excess energy contained in the interface, given that there are no other driving forces preventing this.
Hence in this case we will use \cref{eq:phieq} alone without concentration coupling, which means the sole driving force is the interfacial energy reduction.
Given that curvature flow is the only active mechanism, one may write for the volume of the embedded grain\cite{Balluffi2005}
\begin{align}
  \pdiff{V}{t} &= - \int_{\partial V} v dA\\
  &= - \int_{\partial V} M\gamma\kappa dA
\end{align}
with the grain mobility $M$, the interface energy $\gamma$ and the curvature $\kappa$.
For a particular dimension and object shape the ODE can be solved analytically, say in two dimensions for a circle, resulting in:
\begin{align}
 r(t)^2 - r(0)^2 &= -2M\gamma t\\
 A(t) - A(0) &= -2M\pi\gamma t
\end{align}
which gives us an exact solution to compare to, with the last line not even requiring any shape assumption anymore.
As the time derivative of the area should be constant, the error may naturally be defined as
\begin{align}
 e &= |\pdiff{A_{num}}{t} - (-2M\gamma\pi)|
\end{align}
with the numerical derivative being the average for $t>0$ of a 2nd order finite difference derivative of the respective data.
This time filter is applied to exclude the early adjustment from the ideal to a numerical profile.

The inner grain uses the distance field \cref{eq:centersphere}, with the outer one being the complement.
The area of the grain follows from simple volume integration.
The simulation is terminated upon reaching $t=\num{23e3}$ before the grain vanishes and the curvature becomes comparable to the interface width.

\begin{figure}[h]
  \centering
  \begin{subfigure}[]{0.45\textwidth}
    \includegraphics[width=\textwidth]{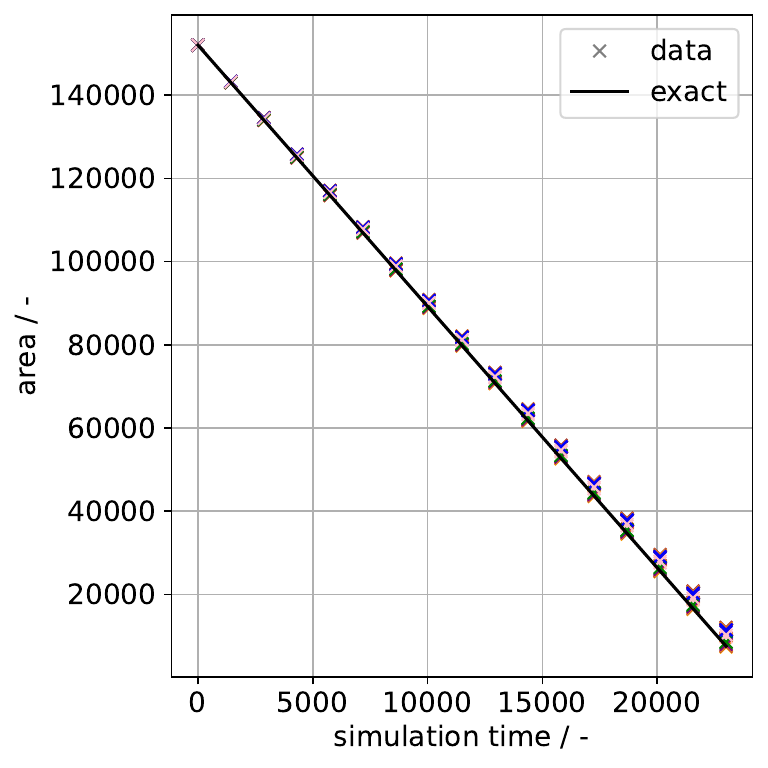}
    \caption{area evolution}
    \label{fig:graingrain-plot}
  \end{subfigure}
  \begin{subfigure}[]{0.45\textwidth}
    \includegraphics[width=\textwidth]{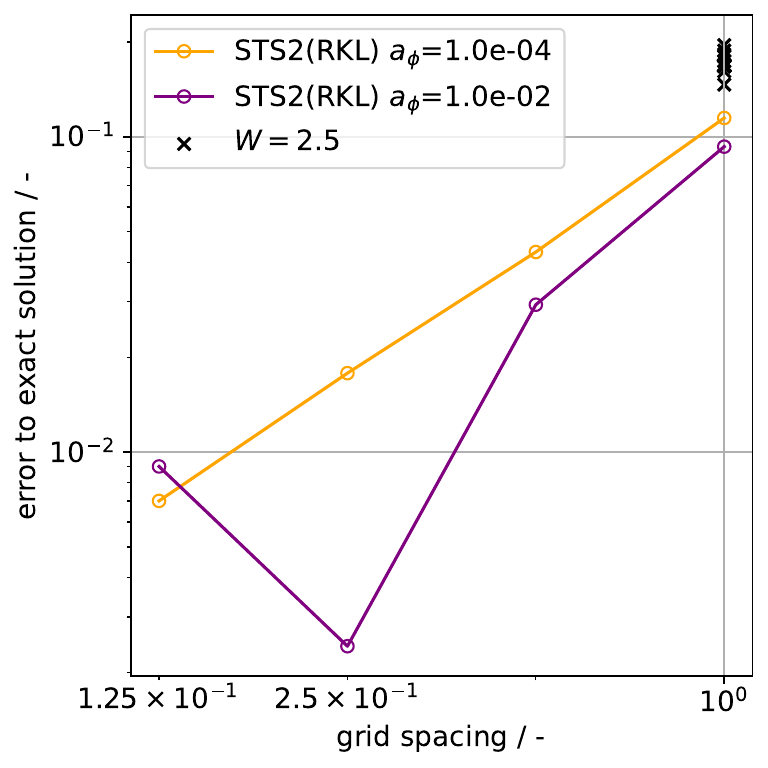}
    \caption{spatial convergence}
    \label{fig:graingrain-errs}
  \end{subfigure}

  \begin{subfigure}[]{0.45\textwidth}
    \includegraphics[width=\textwidth]{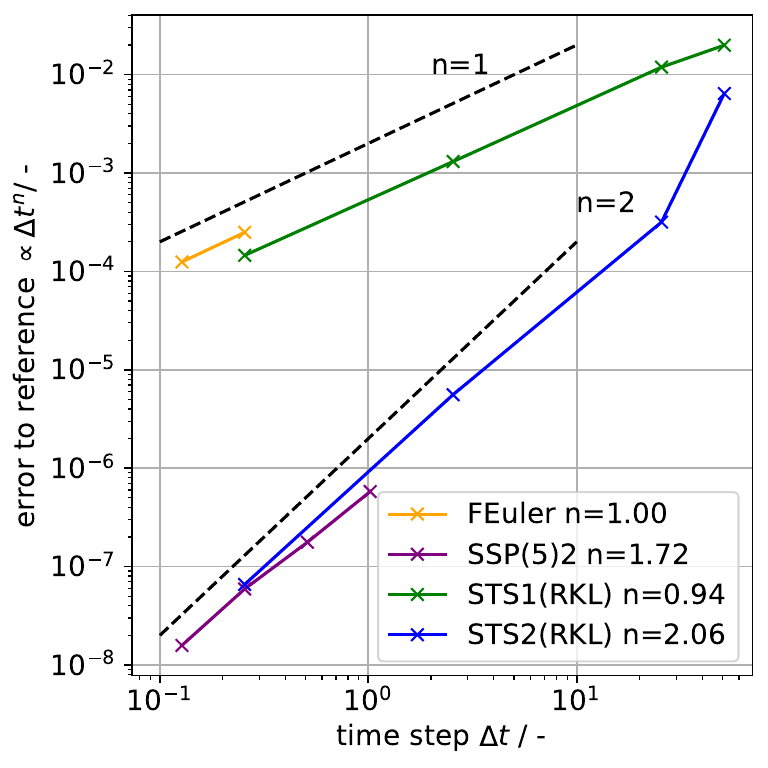}
    \caption{temporal convergence for $W=2.5$}
    \label{fig:graingrain-perf-refsol}
  \end{subfigure}
  \begin{subfigure}[]{0.45\textwidth}
    \includegraphics[width=\textwidth]{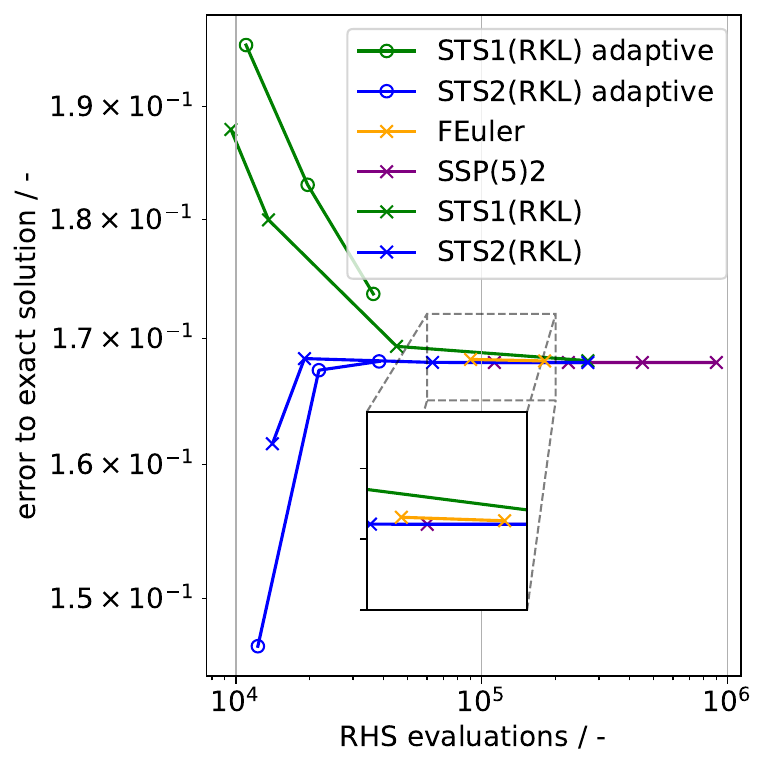}
    \caption{performance for $W=2.5$}
    \label{fig:graingrain-perf-exsol}
  \end{subfigure}
  \caption{Results for the embedded grain geometry. The exact solution is closely approximated by all simulations. The convergence plot shows that for $a_\phi=\num{1e-2}$ the convergence stalls starting between $\Delta x = \{0.25, 0.125\}$, whereas consistent convergence is obtained for $a_\phi=\num{1e-4}$. Comparison to a reference solution shows the expected temporal orders. Comparing to the exact solution shows the STS integrators as the most efficient ones.}
  \label{fig:graingrain}
\end{figure}

A domain of size $512^2$ is employed with an initial seed size of $r=220$ for the embedded grain.
The results of applying the integrators specified in \cref{tab:ints} together with the resolution study are collected in \cref{fig:graingrain}.
In \cref{fig:graingrain-plot}, we can observe that the qualitative time dependence of the exact solution is obtained by all simulations, but with quantitative differences.
In \cref{fig:graingrain-errs} convergence is observed throughout for time integration tolerances of \num{1e-4}, but for $a_\phi = \num{1e-2}$ this eventually limits the convergence.
Interestingly, the error is smaller for the relaxed tolerance prior to this, suggesting that lower error may be achievable using a coarser tolerance due to cancellation of spatial and temporal errors.
Thus for practical simulations, using a relatively coarse absolute tolerance can be beneficial for both speed and accuracy.
Outside of the convergence part, the effect of time integration error is generally small as evidenced by the close clustering of crosses representing the remaining simulations.
The largest of these, corresponding to a STS1 scheme with the largest timestep, has a relative error of about $3\%$, with most other simulations clustering around $2.7\%$ relative error.

Next, looking at only the practical interface width results, we can compare the errors from fixed time step integrators to the reference integrator, which is shown in \cref{fig:graingrain-perf-refsol}.
From this we can observe that generally the expected temporal error orders are obtained\footnote{When used as ODE solvers for the Dahlquist test equation all obtain their order within a tolerance \num{5e-2}.}.
The SSP(5)2 integrator tends to deviate from its order close to its stability limit.
At larger time steps for the STS schemes, the interaction with the exact time matching output tends to force a mix of time steps which can affect the reported order.
By instead comparing to the exact solution we can gauge which is the most efficient scheme:
Plotting the error over the number of RHS evaluations as in \cref{fig:graingrain-perf-exsol} we can observe that generally the STS schemes are the most efficient for this benchmark.
The SSP(5)2 scheme at its stability limit is slightly more precise than the Euler scheme, but also requires slightly more function evaluations, and hence little is gained.
Relative to the FEuler scheme the fixed STS schemes provide up to a factor of 9.4 speedup, whereas the adaptive ones are faster by a factor of 7-8, at little extra error or even smaller error at a slightly smaller speedup.
Fixed timestepping for the STS schemes may outperform the adaptive one here, which is mainly due to a relatively large ratio of rejected steps (up to 29\%).
Note however that with fixed timestepping changing dynamics will not be adapted to, which can lead to inacceptably large errors.

For grain growth with a single interface, the reported speedups should be universal, as the reduced mobility $\gamma M_{\alpha\beta}$ controls the time scale of the problem and the maximal eigenvalue, and hence changing it is only a scaling of time.
Once there are multiple interfaces, as later shown for grain growth, we can expect the speedup for adaptive integrators to decrease, as the size of regions with nontrivial error grows and hence smaller steps must be taken; though at the same time the error committed by the FEuler integrator without error control will grow.
If these interfaces have different mobilities, the speedup would increase however as the increased stiffness can be handled by the STS integrator at lower cost than by the FEuler integrator.

\subsubsection{Solutal Stefan problem}
\label{sec:stefan}
The second kinetic test approximates the classical Stefan problem in a solutal setting, and hence both \cref{eq:phieq} and \cref{eq:ceq} will be used.
Two semi-infinite phases $\alpha$ and $\beta$ are brought into contact with respective initial compositions $c_\alpha$, $c_\beta$.
If these differ from their equilibrium composition, then the interface between them will move.
Following \cite{Balluffi2005}, employing a similarity solution for the diffusion fields and incorporating the flux balance at the interface as well as assuming equal atomic volumes and diffusivities $D$, yields
\begin{align}
 X(t) &= A\sqrt{t} \label{eq:stefangrowth}\\
 A &= \sqrt{\frac{4D}{\pi}} \big{[} \frac{c_\alpha - c_{\alpha\beta}}{c_{\beta\alpha}-c_{\alpha\beta}} \frac{\exp(-A^2/(4D))}{1-\erf(A/\sqrt{4D})}\nonumber \\
  &+ \textcolor{white}{\sqrt{\frac{D}{\pi}} \big{[}} \frac{c_\beta - c_{\beta\alpha}}{c_{\beta\alpha}-c_{\alpha\beta}} \frac{\exp(-A^2/(4D))}{1+\erf(A/\sqrt{4D})} \big{]}
\end{align}
for the interface position $X$ and the growth parameter $A$, with the concentrations on the $\alpha$ and $\beta$ side of the interface $c_{\alpha\beta}$ and $c_{\beta\alpha}$.
Assuming that these interface concentrations match the equilibrium concentrations is equivalent to saying that the generalized Gibbs-Thomson effect vanishes, i.e.
\begin{align}
 \Delta \mu = \Gamma\kappa + K v = 0
\end{align}
with suitably defined Gibbs-Thomson coefficient $\Gamma$ and kinetic coefficient $K$.
Employing a one-dimensional geometry ensures $\kappa$ vanishes.
Since $v>0$ is required for kinetics, we are left with making $K$ vanish by determining its relation to the phase-field parameters.
The latter is usually done with a thin-interface analysis; for details, refer to e.g. \cite{choudhury2012}.
A relationship for the phase-field mobility $L$
\begin{align}
 L &= \frac{\pi^2}{16W^2} \frac{D \pdiff{c}{\mu}}{(c_\beta^{eq} -c_\alpha^{eq})^2 (M+F)} \label{eq:thin-interface}
\end{align}
is eventually obtained for $K=0$, with the terms $M+F \approx 0.3084251$ being integrals involving the weighting function $h$ as well as the phase-field profile.
Note that with the assumption of equal diffusivities there is no need to include an anti-trapping current.
Furthermore, the prefactors $k$ in the Gibbs energies in \cref{eq:gibbs-energy} are reduced to 1 as to allow sufficient solubility in the phases, which is required by the long-ranged smooth concentration fields.

The error is defined with the help of $A^{*}$, obtained from fitting the simulated interface position-time data to \cref{eq:stefangrowth}, as
\begin{align}
 e = |A - A^{*}|.
\end{align}

The distance field for phase $\beta$ is given by $d = h-x$ for an initial height $h$ from the left boundary, with that of $\alpha$ being the complement.
The initial concentration is interpolated as $c(t=0) = \phi_\alpha c_\alpha + \phi_\beta c_\beta$, i.e. the far-field phase concentration values are employed.
For the phase-field, gradient zero boundary conditions are employed, whereas the boundary concentrations are held at the respective initial concentrations with Dirichlet conditions.
The one-dimensional domain's left boundary is held at the equilibrium concentration of $\beta, c_\beta = 0.98$, with the right boundary held at $c_\alpha = 0.2$.
Together with the diffusivity $D=1$ this results in a growth constant $A\approx 0.2412$ via \cref{eq:stefangrowth}.
The initial height $h$ of the $\beta$ phase is 400.
The momentary position $h+X(T)$ is determined based on the position of the $\phi_\alpha = 0.5$ point determined via solving a cubic interpolant.

We note that there is an additional error source here compared to the other benchmark geometries, namely the far-field concentration.
For the analytical solution an infinite domain is assumed, which obviously cannot be simulated.
Hence one must ensure that the gradients in the far field are negligible.
We chose a time giving a certain displacement $X(t)=50 \to t\approx \num{43e3}$ via \cref{eq:stefangrowth} and then varied the physical domain size in an informal study to fix it for the reported simulations, with this ending up being a domain of size 1800.

\begin{figure}[h]
  \centering
  \begin{subfigure}[t]{0.45\textwidth}
    \includegraphics[width=\textwidth]{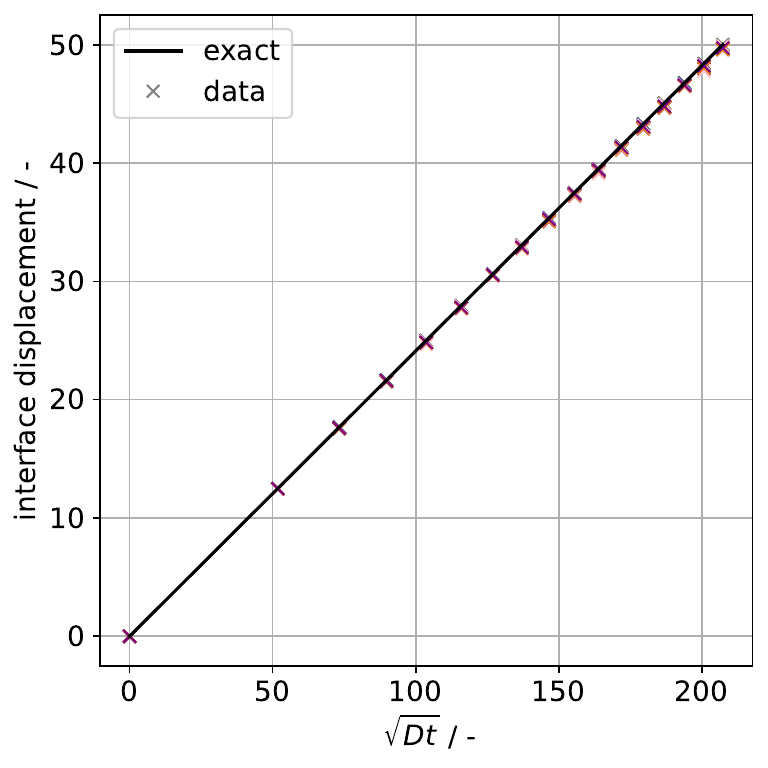}
    \caption{interface motion}
    \label{fig:stefan-plot}
  \end{subfigure}
  \begin{subfigure}[t]{0.45\textwidth}
    \includegraphics[width=\textwidth]{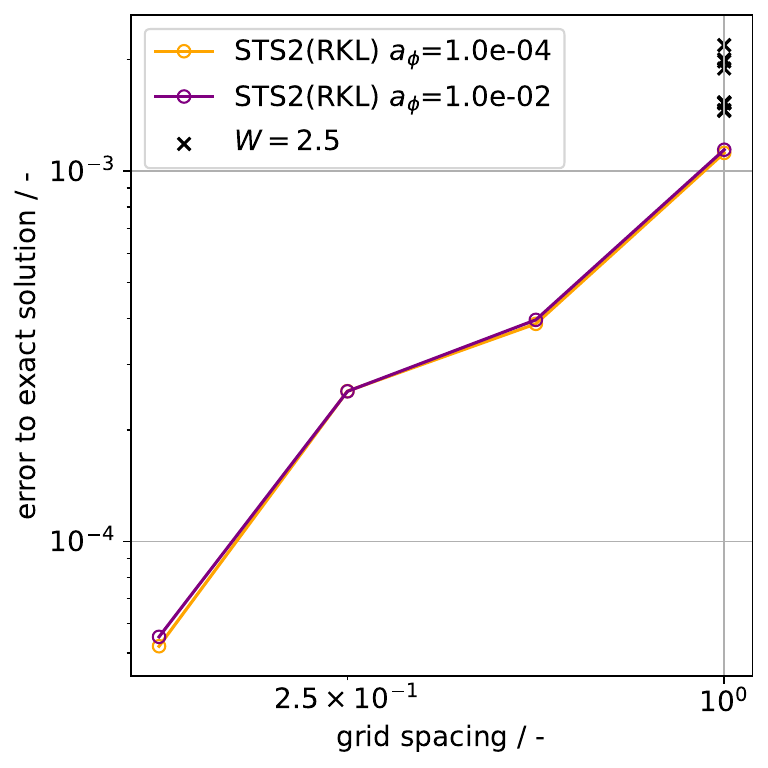}
    \caption{convergence behaviour}
    \label{fig:stefan-errs}
  \end{subfigure}

  \begin{subfigure}[]{0.9\textwidth}
    \includegraphics[width=\textwidth]{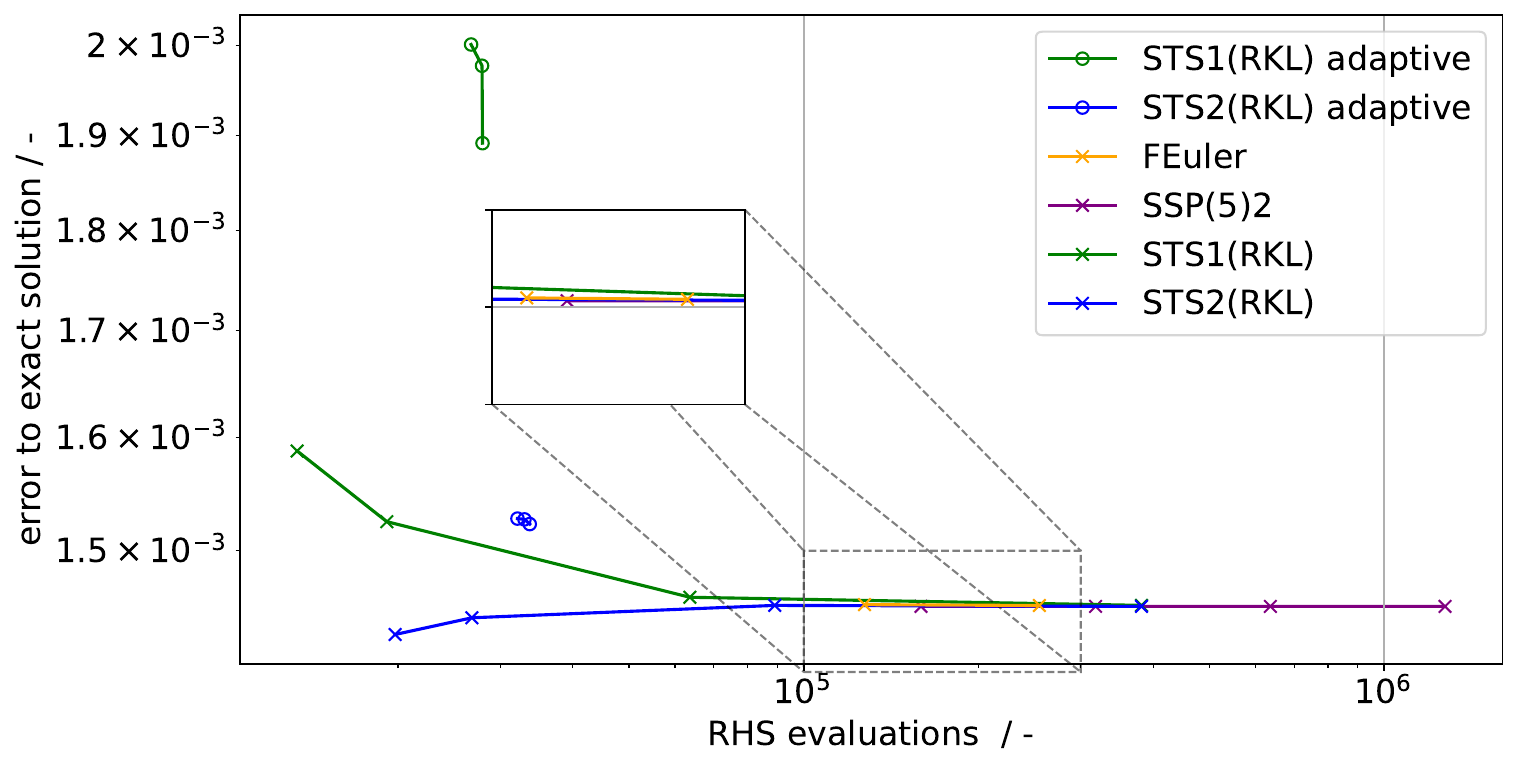}
    \caption{performance behaviour for $W=2.5$}
    \label{fig:stefan-perf}
  \end{subfigure}
  \caption{Results for the Stefan problem. The analytical solution is closely approximated and is being converged to with grid refinement. The choice of integrator seems to have little effect on the error and the STS integrators are observed to be the most efficient ones.}
  \label{fig:stefan}
\end{figure}

All results are collected in \cref{fig:stefan}, in which \cref{fig:stefan-plot} shows a close match to the exact solution.
We observe that the model converges upon grid refinement from \cref{fig:stefan-errs}, with the practical resolution results still showing good accuracy.
All results appear insensitive to the choice of $a_\phi$; informal variation of the remaining tolerances also showed little influence, suggesting the error is mainly controlled by diffuse interface effects, which naturally get scaled down due to the choice of $W(\Delta x)$.
Considering that the rate of change of all field variables is much slower than in the previous case, the lack of sensitivity of the time integration error is quite sensible.

For the practical interface width the work-precision chart is shown in \cref{fig:stefan-perf}.
We observe, as would be suggested by previous results, that the obtained error is largely independent of the integrator and tolerance choice.
Using either of the adaptive STS schemes here gives a speedup of about 4-5 with effectively no additional error; fixed timestepping produces speedups from 6 to 10.
The phase-field tolerance seems to play little role for the performance here, suggesting that the temporal error is dominated by the concentration tolerance which was fixed for this study.
Setting $a_c = a_\phi$ showed a performance improvement (up to a factor of 11) at no significant change in error, as shown in the supplementary material.

The speedups reported here depend on the Gibbs free energies by relating to the mobility \cref{eq:thin-interface} but also on the diffusivities if these are allowed to vary between the phases, though in this case a quantitative phase-field simulation requires the addition of some model term counteracting the interface width dependent excess solute trapping.
If this term introduces extra stiffness, the resulting speedups would again be larger, given that the time scale of the problem and hence the magnitude of its time derivatives are controlled by the sharp interface problem independent of the phase-field model.
Furthermore, a larger diffusional contrast is beneficial for the speedup as well, since the error in the low-diffusion regions will be smaller and hence an overall larger timestep may be used to achieve the same error.

Finally, none of the kinetic benchmarks showed any energy increase.
It is likely that any actual energy error is drowned out by the magnitude of energy dissipation for these cases.

\subsection{Grain growth with particles}
\label{sec:ggpores}
As an application example we will also treat a system which doesn't have an analytical solution, but which is representative of more typical ratios of kinetic parameters as well as more complicated geometries.
For this a system of many grains $\beta$ of different orientation is employed, together with spheres of a second phase $\alpha$ randomly placed in this structure.
When these second phases are in contact with a grain boundary, they exert a drag force on the boundary.
Depending on the effective mobility of the second phase, these may be taken to be immobile particles which act as pinning centers, allowing for investigations relating to Zener pinning and limiting grain sizes.
For more mobile phases, this represents the final stage of sintering with only closed porosity remaining, which reduces grain growth but does not necessarily lead to a limiting grain size, as pores may move and end up merging, becoming less effective pinning centers.
The effective particle mobility can be controlled via the size of the second phase and the diffusivity controlling the mass transport\cite{Liu1994}, which for simplicity here was taken to be diffusion within the second phase, i.e. vapour diffusion if we interpret it as a pore.
We note that Zener pinning can in general be more easily simulated without recourse to using a concentration field, see e.g. \cite{Harun2006}.
The final sintering stage can also be approximated without diffusion fields, see e.g. \cite{Rehn2019}, but this assumes known dominant diffusion mechanisms, more materials parameters and a tracking algorithm to account for pore pressure and pore-pore interactions.
Furthermore, pores in the final sintering stage can also be elongated, in which case diffusional instabilities such as the Plateau-Rayleigh instability become relevant.
These are more easily described by employing concentration fields.
We note that the eventually chosen parametrization for this case may be calculated more efficiently by integrating the phase-field and concentration fields separately, as their respective timescales differ by several orders of magnitude; but for ease of presentation and comparability both fields were integrated together.

We seek to reproduce the following qualitative behaviours with three simulations:
If grain growth without particles happens, mean-field theory based on curvature as the sole driving force predicts parabolic growth of the mean grain size, i.e. $G(t)^n - G(0)^n = kt, n=2$\cite{Rahaman2017}.
As alluded to in the previous paragraph, immobile particles can lead to a limiting grain size, which in experiments may also be identified as large exponents $n$ since the curve will tend to flatten out increasingly.
In contrast to this, while mobile particles (pores) also exert a drag force on the grain boundary, they co-move with the grain boundary and may merge with others, reducing their effectiveness as pinning centers.
Hence in this case the exponent generally increases, see e.g. \cite{Nichols1966,Riedel1993}, while not necessarily reaching a limiting grain size.
Based on the different parameters here we will be able to observe speedups for realistic parameter combinations.

The initial conditions for pure grain growth are a periodic Voronoi tesselation of 3D space with random points, with a minimal distance of $12\Delta x$ between each point of a $(400\Delta x)^3$ domain.
Once this simulation has reached a roughly parabolic regime ($t=1$ in the computational experiments) simulations with particles are constructed based on this configuration:
For immobile particles, spheres are placed randomly in space with a target volume ratio $f_i=0.15$ and radius $r_i=6\Delta x$ in order to showcase the limiting grain size, with the distance between spheres being at least $d=2r$.
Due to the finite extent of the phase-field this still means that some overlap exists, which does increase the effective radius but due to the low mobility the particles do not move or change shape over the simulation duration.
However, for mobile particles this would lead to quick merging of most particles, giving a completely different geometry than expected.
Hence for these a minimal distance of $d=2(r+e), e=8\Delta x$ is imposed, which limits the obtainable porosity as a function of this distance and the pore radius.
Due to this limit we pick $f_m=0.05, r_m = 10\Delta x$ which is also more representative of the porosity in final stage sintering.
While it is possible to also calculate the immobile particle case with this particle distribution, the simulation domain would need to be enlarged significantly to reach the limiting grain size.
For the immobile parameters the classical Zener pinning theory predicts a limiting grain size of $r_{l,i} = \frac{4}{3}\frac{r}{f} \approx 53.3\Delta x$ whereas for the mobile particle distribution it would be $r_{l,m} \approx 266\Delta x$ and hence require a domain in excess of $(800\Delta x)^3$\cite{Miodownik2000}.
While calculating such a domain poses no problems but the waiting time for the employed code, the goal of reproducing qualitative behaviour does not necessarily justify the simulative expense.

We note that the phase-field volume $\int_V \phi(d)dV$ of a two and higher dimensional shape has an $O(W^2)$ error compared to the sharp interface volume the distance field $d$ is supposed to represent.
In order to still obtain a porosity close to the target porosity even when the size of the particles is comparable to the interface width, the radius $r$ is transformed to a smaller equivalent radius $r_c$ such that the target volume is obtained with the help of a function function $r_c=f(r,W)$ derived in the supplementary material.
With the centers and radius of the particles fixed, the only thing left is to place them in the domain.
Since a diffuse initial shape is used, $\phi_\alpha = \phi(r_c-x_{min})$ is imposed for each cell while using \cref{eq:distfix} to update the already present phase-fields.
$x_{min}$ here is the smallest distance of this cell to any of the particle positions.
The geometry is hence fully described and what remains is to fix the parameters representing the different regimes.

For insoluble, mobile particles (pores) we want small particles in some sense with sufficient diffusivity.
Since the mobility of these scales as $r^{-n}, n>2$ \cite{Liu1994} with the particle size $r$, choosing a smaller length scale means we have more mobile particles relative to pure grain growth.
At the same time, we would like to use the same length scale for immobile particles, hence after fixing the length scale, diffusion values leading to immobile and mobile particles are determined.
Furthermore, in order for the particles to be insoluble even under capillary action the prefactor $k$ in \cref{eq:gibbs-energy} needs to be large relative to capillary forces.
We can estimate this effect by considering the expression for the phase concentration in terms of the chemical potential $c_\alpha = c_{0,\alpha} + \frac{\mu}{k}$ and the equilibrium condition $\psi_\alpha - \psi_\beta = \gamma\kappa$, from which we find that the shifted equilibrium bulk concentration is $c_\alpha(\gamma\kappa) = c_{0,\alpha} - \frac{\gamma\kappa}{ k(c_{0,\alpha} -c_{0,\beta}) }$.
We would like for the difference to the bulk concentration to be small relative to the bulk concentration, and hence must know the scale of $O(\kappa) = O(\frac{1}{r})$.
By computational experiments we pick a length scale of $\Delta x = 0.005$ at which $D_\alpha=10$ shows visibly mobile particles and hence $k = \num{1e5}$ ensures $\frac{c_\alpha(\gamma\kappa) - c_{0,\alpha}}{c_{0,\alpha}} < 1\%$ for the considered particle sizes.

The grain boundary mobility is set to be $M_{\beta\beta} = \num{5e-3}$.
The diffusivity within grains is taken to be $D_\beta = \num{5e-8}$ as a small number effectively suppressing Ostwald ripening on the considered time scale.
For immobile particles $M_{\alpha\beta} = \num{5e-8}, D_\alpha=\num{5e-8}$ are chosen as to make the particles immobile on grain growth time scales.
Finally, with $M_{\alpha\beta} = \num{5e-3}, D_\alpha=10$ mobile particles are realized.

If the immobile particle parameters are employed directly after setting the particles, stable many phase regions tend to form on the fringes of the particles, even if a simplex accounting for mobility variations is employed as suggested in \cite{Daubner2023}.
We suspect that this is due to chemical driving forces aiming to shrink the particles to achieve their capillary equilibrium, which gives a driving force to any grain phases in the vicinity.
This driving force may exceed the potential energy terms of the phase-field and hence stabilize many phase regions which originally exist due to the interface width being comparable to some of the initial grains' equivalent radius.
While these would eventually vanish after capillary equilibrium is achieved, the timescale for that is far beyond the timescale of grain growth for this parameter combination.
Since these may act as artificial pinning centers, an additional equilibration run is conducted for both particle-laden simulations prior to running with the actual parameters.
The run uses $M_{\beta\beta} = \num{5e-8}, D_\alpha=\num{5e-3}, M_{\alpha\beta}=\num{5e-3}$ as to ensure that virtually no grain growth happens while at same time allowing for capillary relaxation, and is continued for a duration of $5$ time units.
After this equilibration run the time is reset to $t=1$ and the simulations are run normally; no further formation of stable many phases regions was observed.
All runs used the tolerances $a_\phi=\num{1e-2}, a_c = \num{1e-3}$, with the latter being chosen as early tests showed it had little influence as in the Stefan benchmark.

\begin{figure}[h]
    \includegraphics[width=\textwidth]{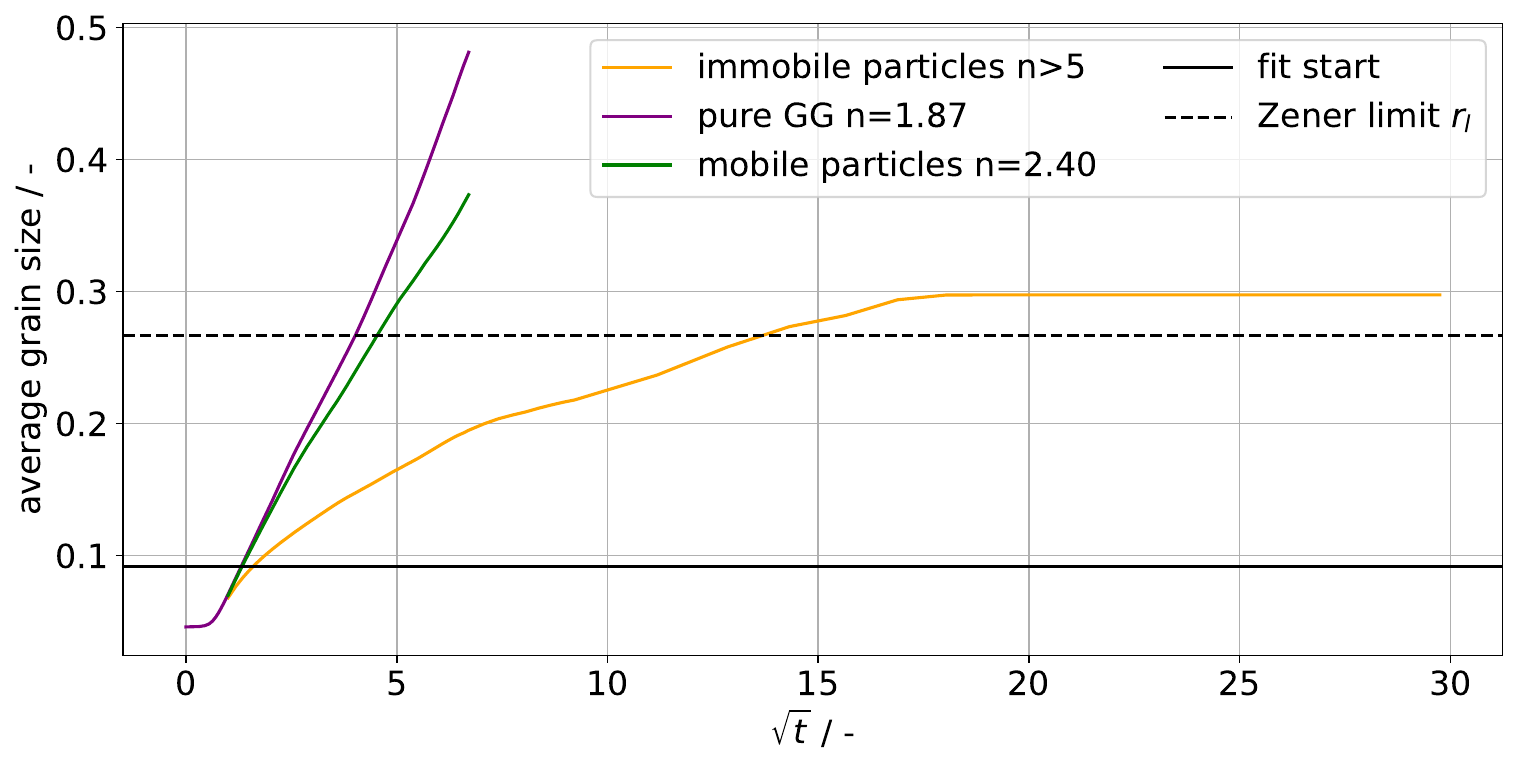}
  \caption{Results for grain growth, either pure, with immobile particles or mobile particles. The expected qualitative behaviour is observed, with exponents close to expectations being obtainable. For immobile particles, a limiting grain size is eventually reached.}
  \label{fig:ggcurves}
\end{figure}

The three simulations, conducted as described above, yield the grain size curves in \cref{fig:ggcurves}.
An exponent close to $n=2$ is observed for pure grain growth, whereas an increased exponent is observed for grain growth with mobile particles, as would be expected.
Theoretically one would expect $n=3$ for mass transport through the vapor with the pore pressure being maintained at capillary equilibrium $p=\frac{2\gamma}{r}$\cite{Nichols1966}.
We suspect that this discrepancy is due transient behaviour between the pure and pore-dragged behaviour, as the initial part of the curves are quite similar; reaching $n=3$ likely requires significantly larger domains to allow for the dynamics to emerge.
For the immobile particles we can compare the obtained limiting grain size (constant for about a time span of 200) to the prediction of Zener, which comes out to be $r_l\approx0.27$.
The simulation predicts a larger value of about $0.3$, which is likely due to the interface width being of the same order of magnitude as the radius of curvature and overlap between particles effectively increasing the average radius.
Hence we can say that even at quite coarse tolerances the STS schemes still recovers the essential features of the free-boundary problem we are concerned with.

Selected simulation snapshots are shown in \cref{fig:gg3d-imgs} with videos of the evolution deposited in the supplementary material.
The immobile particle case is depicted with only the grain boundary (GB) network as the 0.5 isolevel of the $\phi_{gb} = \min(4\sum_{i\neq\alpha, j\neq\alpha,j>i} \phi_j\phi_i, 1)$ field.
The GBs are disturbed by circular holes due to the immobile particles and hence the particles' position may be inferred from these; note that the radius of the holes depends on both the particle radius and where the GB intersects them.
Once a GB has accumulated sufficient particles, it apparently stops as highlighted; other faces of the grain may still advance however.
The limiting grain size is reached once all GBs have a sufficient particle density.

\begin{figure}
  \centering
  \includegraphics[width=0.9\textwidth]{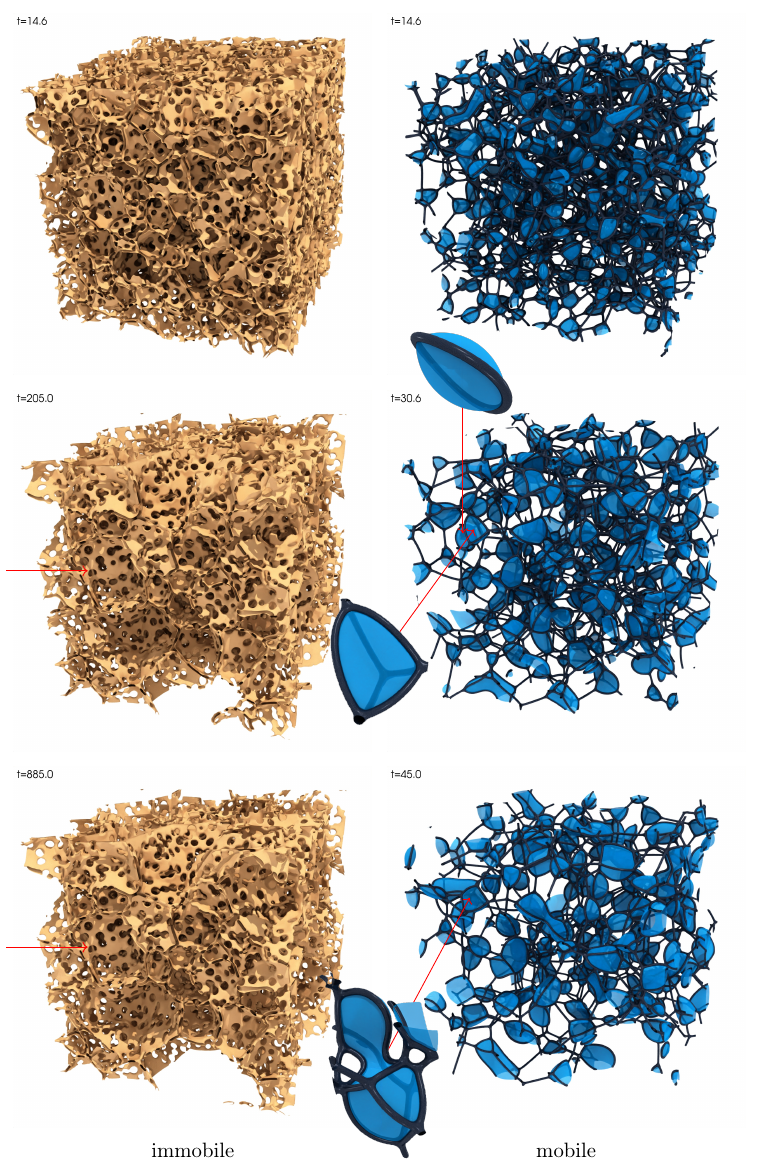}
  \caption{Selected snapshots of the 3D evolution with particles. For immobile particles (left column) GBs eventual stop their motion once they are sufficiently laden with particles as indicated. For mobile particles (pores) on the right column, these move and merge over time, allowing for continued, albeit slowed down, grain growth. Pores attain a variety of shapes depending on their surroundings and time to equilibrate.}
  \label{fig:gg3d-imgs}
\end{figure}

The simulation with mobile particles show the triple lines (dark) defined as the 0.5 isolevel of the $\phi_{tj} = \min(27\sum_{i,j>i,k>j}\phi_i\phi_j\phi_k, 1)$ field together with the 0.5 isolevel of the particle phase (bluish transparent).
The pores co-move together with the grain boundaries and merge over time, reducing their number.
As indicated by the closeup views, pores can attain a variety of shapes depending on their environment and time to equilibrate.
If the pores only touch two grains, their shape tends to be a squashed sphere, though most pores touch more grains.
In that case the shape is close to a rounded polyhedron with vertices located at the enclosing grain triple lines in order to satisfy the surface tension balance; a tetrahedral example is shown for the middle row.
This is in line with \cite{Tikare2001d,HOTZER2016}, who employed a more regular geometry with pores being placed on quadruple points between a basal hexagonal polycrystal and one top crystal growing into this polycrystal.
They observed that the pores took on shapes similar to ice-cream cones, with a triangular cross-section, which is the cross-section of a tetrahedron.
Pores merging results in temporary irregular shapes as shown on the bottom.
This pore is the result of the tetrahedral pore (four-sided face showing) from the middle row merging with the spheroidal pore shown in the middle panel and one nearby other tetrahedral pore.
A video of the time evolution of this subdomain is deposited with the supplementary material.
Due to the intentionally large pore mobility no pores are left behind within the growing grains.

The performance results are shown in \cref{tab:perf-ggpores}.
The RHS count for the FEuler integrator is estimated as $\frac{t_{r}}{\Delta t_{e}}$ as, like the table shows, the mobile pore case would be infeasible to compute with it\footnote{The STS2 simulation took about 16 hours on 4 NVIDIA A100 GPUs; the equivalent FEuler simulation would have taken about 80 days. Whether the system would have run uninterrupted for that long is another question.}.
The duration $t_r$ is the end time minus any time offset, which for both particle-laden simulations is an offset of $1$.
For pure grain growth (GG), $\Delta t_{e}=\num{8.44e-4}$ is the stable Euler time step as the chemical driving force related terms vanish.
For immobile particles this reduces to $\Delta t_{e}=\num{6.26e-5}$ via these terms.
Finally, for mobile particles the stable step is $\Delta t_{e}=\num{2.82e-7}$ via the diffusion equation.
For pure grain growth, the speedup is somewhat smaller as many phases dying off tends to force smaller time steps to accurately resolve the transition.
However, once the problem becomes stiffer, this is more than compensated by the increased stability range of the integrator.
\begin{table}[h]
\centering
\caption{3D grain growth performance in terms of RHS evaluations}
\label{tab:perf-ggpores}
 \begin{tabular}{l|ccc}
 & pure GG & immobile particles & mobile particles \\
\hline\\
duration $t_r$ & 45 & 884 & 44\\
FEuler & \num{53279} & \num{14110453} & \num{156288000}\\
STS2 & \num{14288} & \num{1224604} & \num{1365954} \\
speedup & \num{3.73} & \num{11.5} & \num{114} \\
\end{tabular}

\end{table}

For pure grain growth, the speedup can likely be seen as a lower limit for the simplest phase-field problems, as the large amount of evolving interfaces requires fine temporal resolution to solve accurately before transitioning into a more regular behaviour.
Adding mobility or interfacial energy contrasts or anisotropy will increase the speedup as the stiffness increases.
External coupling to a known temperature ramp to reproduce experimental results would achieve even larger speedups since while the external timescale is fixed, the STS integrator can take dynamics-appropriate timesteps to cover the external timescale.\\
With particles, the diffusion value enters the stiffness question, but at the same time the error in the concentration evolution will be affected too.
Seeing the immobile particles case as the zero diffusion limit of the mobile particles case, this gives us a lower limit for the speedup, with increasingly larger diffusion coefficients increasing the speedup.

\FloatBarrier
 
\section{Conclusion \& outlook}
Within the work it could be shown that energy stability has little influence on the error w.r.t. the sharp interface solution of two benchmark problems.
Based on this, kinetic verifications showed behaviour convergent to sharp interface solutions, with practical choices of interface width still closely approximating these.
Overall, both spatial and interfacial width related errors tended to dominate even at coarse integration tolerances; the phase-field could generally be integrated with absolute tolerances up to $\num{1e-2}$ without significant accuracy loss at practical spatial resolution.
With adaptive integrators, the temporal error only became visible after several levels of grid refinement.
Speedups obtained range from a factor of 4 to 114 when using a STS integrator over the forward Euler integrator.
This becomes the most pronounced when the equations exhibit multiple time scales, as is common in many manufacturing processes.
The SSP(5)2 integrator --- while integrating the time derivative more accurately --- does not offer a speedup relative to the forward Euler integrator, as the spatial error dominates over the temporal error for practical resolutions.
This likely extends to other ``normal'' explicit Runge-Kutta integrators as their stability domain scales linearly in the number of stages at best.
As the STS integrator is not much more complicated than the forward Euler integrator to implement, there is little reason even for self-rolled codes not to prefer it over the forward Euler integrator.

If advective terms are included, the speedups would likely even be larger as accurate advection schemes such as WENO are relatively expensive to compute.
By either using a splitting scheme or a partitioned scheme, these can be integrated to solve the evolution equations with the appropriate stability domains, while reducing the number of required number of advection term evaluations due to requiring fewer steps.
The same idea applies to coupling to stationary problems e.g. stationary mechanics, which will need to be solved for for far fewer times.

The speedups of the STS schemes relative to the FEuler scheme given in the work can typically be seen as lower limits for the respective benchmark cases and the applications which they represent.
The benchmarks relating to grain growth can be taken as precursors to sintering, where advective terms need to be added as described above.
The Stefan problem is the prototype for phase-changes problems such as dendritic and eutectic solidification in higher dimensions; the speedup reported here should readily transfer to these.
As the time scale on which dynamics occurred was not far removed from the time scale induced by the largest eigenvalue of the right-hand side, there was no extra stiffness biasing the results against the FEuler scheme.
Only for the case of mobile pores did these time scales differ substantially, which resulted in larger speedups.

What bears further investigation is the relatively large number of step rejections for some benchmarks.
This suggests that either the stability limit is passed or that some part (error estimator, time step adapter) in the step adaption machinery fails.
The former might be probed by including a dynamical estimate of the largest eigenvalue of the system.
If this doesn't affect the rejection rate significantly, it might be that standard adapters need significant tuning in their coefficients to reduce their rejection rate as shown by \cite{Ranocha2022} for compressible CFD problems.

A follow-up paper investigating both explicit and implicit schemes with a well potential is planned to complement the present limited investigation.
 
\section*{Author contributions}
\textbf{Marco Seiz}: Conceptualization, Software, Methodology, Investigation, Data Curation, Validation, Visualization, Writing - original draft, Writing - review \& editing.
\textbf{Tomohiro Takaki}: Funding acquisition, Project administration, Resources, Writing - review \& editing.

\section*{Conflicts of interest or competing interests}
The authors declare that there are no conflicts of interest.

\section*{Data and code availability}
The raw data is deposited with the supplementary material, except for the 3D simulation data for which only the processed data is deposited.
It also contains electronic notebooks producing plots and auxiliary information.
The 3D data may be shared upon reasonable request.
The code required to reproduce the present work cannot be shared publicly.

\section*{Supplementary Material}
The supplementary material of this paper is available at \url{https://zenodo.org/records/18320716}.

\section*{Acknowledgements}
This paper is based on results obtained from a project, JPNP22005, commissioned by the New Energy and Industrial Technology Development Organization (NEDO).

\bibliography{literature}

\end{document}